\documentclass[11pt]{article}

\usepackage[english]{babel}
\usepackage[utf8x]{inputenc}
\usepackage[authoryear]{natbib}

\usepackage{booktabs,bm,soul}
\usepackage[flushleft]{threeparttable}
\usepackage[T1]{fontenc}
\usepackage{algorithm}
\usepackage{algpseudocode}
\usepackage{float}
\usepackage{multirow}
\usepackage[letterpaper,margin=1in]{geometry}
\usepackage[hidelinks]{hyperref}
\usepackage{enumerate}

\usepackage{amsmath}
\usepackage{cleveref}
\usepackage{graphicx}
\usepackage{breqn}
\usepackage[dvipsnames]{xcolor}

\usepackage[small,compact,noindentafter]{titlesec}
\titlespacing{\section}{0pt}{2ex}{0ex}
\usepackage[colorinlistoftodos]{todonotes}

\usepackage{amssymb}
\usepackage{amsthm}
\makeatletter
\newcommand{\newreptheorem}[2]{\newtheorem*{rep@#1}{\rep@title}\newenvironment{rep#1}[1]{\def\rep@title{#2 \ref*{##1}}\begin{rep@#1}}{\end{rep@#1}}}
\makeatother

\usepackage{mathrsfs}
\usepackage{blindtext}
\usepackage{cancel}
\usepackage[mdyy,12hr]{datetime}
\usepackage{scalerel} 
\usepackage{rotating,tabularx}
\usepackage{subcaption}
\usepackage{mwe}
\usepackage{url}
\usepackage{bbm}

\usepackage{enumitem}
\setenumerate[0]{label=(\arabic*)}
\usepackage[mathscr]{eucal}
\theoremstyle{remark}
\newtheorem{innercustomgeneric}{\customgenericname}
\providecommand{\customgenericname}{}
\newcommand{\newcustomtheorem}[2]{%
	\newenvironment{#1}[1]
	{%
		\renewcommand\customgenericname{#2}%
		\renewcommand\theinnercustomgeneric{{\it ##1}}%
		\innercustomgeneric
	}
	{\endinnercustomgeneric}
}

\newcustomtheorem{problem}{Problem}
\newtheorem{thm}{Theorem}
\newreptheorem{theorem}{Theorem}
\newtheorem{lemma}{Lemma}
\newreptheorem{lemma}{Lemma}

\newtheorem*{lemma*}{Lemma}
\newcustomtheorem{sublemma}{Sublemma}
\newcustomtheorem{subsublemma}{Subsublemma}

\newreptheorem{proposition}{Proposition}
\newcustomtheorem{theorem}{Theorem}
\newtheorem*{theorem*}{Theorem}
\newtheorem*{remark}{Remark}

\newcommand{\be}{ \begin{equation} }
	\newcommand{\ee}{ \end{equation} }

\DeclareMathOperator{\supp}{supp}

\renewcommand{\phi}{\varphi}

\mathchardef\hyphen="2D 

\algdef{SE}[DOWHILE]{Do}{doWhile}{\algorithmicdo}[1]{\algorithmicwhile\ #1}%

\makeatletter
\def\smallunderbrace#1{\mathop{\vtop{\m@th\ialign{##\crcr
				$\hfil\displaystyle{#1}\hfil$\crcr
				\noalign{\kern3\p@\nointerlineskip}%
				\tiny\upbracefill\crcr\noalign{\kern3\p@}}}}\limits}
\makeatother
\usepackage{setspace}
\allowdisplaybreaks
\begin{document}
	
	{\centering
		{\bf\Large Closing the Gap: Efficient Algorithms for Discrete Wasserstein Barycenters}\vspace{0.1cm}
		\\{Jiaqi Wang}
		\\{\small H. Milton Stewart School of Industrial and Systems Engineering, Georgia Institute of Technology, Atlanta, GA, USA, jwang3737@gatech.edu}
		\\{Weijun Xie}
		\\{\small Department of Data Science, City University of Hong Kong, Hong Kong, wj.xie@cityu.edu.hk}
		\\
	}

	\begin{abstract}
		Given a collection of probability measures, the Wasserstein barycenter problem seeks a probability measure that minimizes the weighted sum of their squared type-$2$ Wasserstein distances. We study the discrete setting, in which every input measure has finite support; this setting arises frequently in machine learning and operations research. Since computing a discrete Wasserstein barycenter is NP-hard, we develop approximation algorithms with provable performance guarantees. The best previously known polynomial-time guarantee is a factor of $2$. We improve this bound by proposing a polynomial-time approximation scheme (PTAS) that generalizes the existing $2$-approximation: for any prescribed $\alpha\in(0,1]$, the proposed algorithms return a $(1+\alpha)$-approximate barycenter in time polynomial in $(nk)^{1/\alpha}$ and $d$. We provide both randomized and deterministic constructions and derive a sharper guarantee when the input measures are equally weighted. Numerical experiments on synthetic and real datasets show that the algorithms are computationally practical and produce near-optimal barycenters.
	\end{abstract}

	\section{Introduction}
	
	A Wasserstein barycenter is a Fr\'echet mean in Wasserstein space: it minimizes a weighted average of optimal-transport costs over a collection of probability measures. For $\mathbb{P},\mathbb{Q} \in \mathcal{P}_2(\mathbb{R}^d)$, the type-$2$ Wasserstein distance is
	\[
	W_2(\mathbb{P},\mathbb{Q})
	:= 
	\left(
	\inf_{\Pi \in \mathcal{M}(\mathbb{P},\mathbb{Q})}
	\int_{\mathbb{R}^d \times \mathbb{R}^d}
	\|\bm{x}-\bm{y}\|_2^2 \, d\Pi(\bm{x},\bm{y})
	\right)^{1/2},
	\]
	where $\mathcal{M}(\mathbb{P},\mathbb{Q})$ denotes the set of couplings with marginals $\mathbb{P}$ and $\mathbb{Q}$.

	Let $[k]:=\{1,\ldots,k\}$, and let $\mathbb{P}_1,\ldots,\mathbb{P}_k$ be discrete probability measures. For each $i\in[k]$, the measure $\mathbb{P}_i$ is supported on a finite set $\Xi_i\subseteq\mathbb{R}^d$ with $|\Xi_i|\leq n$. We denote its $j$-th support point by $\hat{\bm{x}}_{ij}$ and the corresponding mass by $\hat p_{ij}$, for $j\in[|\Xi_i|]$. Given weights $\bm{\lambda}=(\lambda_1,\ldots,\lambda_k)$ in the simplex $\Delta_k:=\{\bm{\lambda}\in\mathbb{R}_+^k:\sum_{i\in[k]}\lambda_i=1\}$, the discrete Wasserstein barycenter problem is
	\begin{equation}\label{eq:def_barycenter}
		v^*=\inf_{\mathbb{P} \in \mathcal{P}(\mathbb{R}^d)}  \sum_{i\in [k]} \lambda_i\, W_2^2  \left(\mathbb{P}, \mathbb{P}_i \right),\tag{WBCenter}
	\end{equation}
	where $\mathcal{P}(\mathbb{R}^d)$ denotes the set of Borel probability measures on $\mathbb{R}^d$ with finite second moment. Discrete Wasserstein barycenters have applications in clustering \citep{ye2017fast}, regression \citep{bonneel2016wasserstein}, dictionary learning \citep{schmitz2018wasserstein}, texture mixing \citep{bonneel2015sliced}, image morphing \citep{simon2020barycenters}, medical imaging \citep{janati2020multi}, facial detection \citep{yan20212d}, time-series modeling \citep{cheng2021dynamical}, and distributionally robust optimization \citep{lau2022wasserstein}.

	\subsection{Literature review}
	The Wasserstein barycenter was introduced by \cite{agueh2011barycenters}, who also connected it to the multi-marginal optimal transport (MOT) problem. Several extensions and variants have since been studied; see, for example, \cite{bigot2012consistent,huang2021projection,uribe2018distributed}. For finitely supported measures, \cite{anderes2016discrete} proved fundamental structural properties--including discreteness, sparsity, and a non-mass-splitting property--and derived an LP formulation through MOT. The resulting formulation, however, has $O(n^k)$ variables and constraints and is therefore exponential in the number of marginals $k$. Hardness persists under substantial restrictions: \cite{borgwardt2021computational} proved that the sparse barycenter problem is NP-hard even for $d=2$ and $k=3$, while \cite{altschuler2022wasserstein} ruled out a polynomial-time randomized algorithm that achieves prescribed additive accuracy, even under uniform weights, unless $\mathrm{NP}\subseteq\mathrm{BPP}$. The latter result concerns additive approximation and thus does not preclude the multiplicative guarantees studied here. Existing multiplicative guarantees are limited to a factor of $2$: \cite{borgwardt2022lp} used the union of the input supports as the candidate barycenter support, and \cite{lindheim2023simple} constructed a barycenter from a sequence of two-marginal transport plans. We improve this factor-$2$ bound by developing a polynomial-time approximation scheme (PTAS): for every $\alpha\in(0,1]$, our method computes a $(1+\alpha)$-approximate barycenter in time polynomial in $(nk)^{1/\alpha}$ and $d$.
	
	Exact algorithms are available under additional structural assumptions or through more elaborate LP techniques. In fixed dimension $d$, \cite{altschuler2021} showed that an exact type-$p$ Wasserstein barycenter for $p\in\{1,2\}$ can be computed in time polynomial in $n$, $k$, and $\log U$, where $\log U$ bounds the input bit length. The same framework computes a barycenter with additive error $\epsilon>0$ in time polynomial in $n$, $k$, and $\log(1/\epsilon)$ by solving an exponential-size LP through a separation oracle based on power diagrams. Outside the fixed-dimensional setting, \cite{borgwardt2020improved} strengthened the standard LP formulation by exploiting the non-mass-splitting structure, aggregating transport variables by support tuples, and incorporating preprocessing rules and support-size bounds. Building on these ideas, \cite{borgwardt2022column} developed an exact column-generation method whose pricing problem searches for barycenter atoms with negative reduced cost and is accelerated through bounding, heuristics, stabilization, and warm starts. These approaches can be effective in practice, but they do not yield polynomial-time algorithms for all general input parameters and can become difficult to scale when $n$, $k$, and $d$ are large.

	Related literature studied the fixed-support variant, in which the barycenter support is specified \emph{a priori} and the remaining problem is an LP. Since this LP can still be large, \cite{cuturi2014fast} introduced entropic regularization and subgradient-based schemes, \cite{benamou2015iterative} developed a parallelizable iterative Bregman-projection (Sinkhorn) method, and \cite{janati2020debiased} proposed a debiased variant that removes the regularization bias while retaining comparable convergence rates. Other approaches include an interior-point method \citep{ge2019interior}, a deterministic Bregman-projection algorithm with improved complexity guarantees \citep{lin2020fixed}, and a symmetric Gauss--Seidel ADMM derived from the dual formulation \citep{yang2021fast}. Although some of these methods can be incorporated into schemes for updating non-prespecified support locations, their 
	analysis does not provide the globally guaranteed candidate-support reduction 
	developed in this paper.

	Several methods address free-support or regularized barycenter problems directly. \cite{claici2018stochastic} proposed a stochastic scheme that updates the barycenter support using sampled subgradients, and \cite{luise2019sinkhorn} developed a Frank--Wolfe method for barycenters defined through the Sinkhorn divergence. \cite{lin2025federated} formulated the free-support problem as an integer program and designed a federated dual-subgradient algorithm, while \cite{kroshnin2019complexity} analyzed iterative Bregman projections and accelerated gradient methods for the entropically regularized objective. These methods provide useful computational tools, but they do not show the polynomial-time multiplicative approximation guarantees developed in this paper.

	\subsection{Summary of contributions}
	We develop a unified candidate-support reduction framework for the discrete Wasserstein barycenter problem \eqref{eq:def_barycenter}. The framework converts the design of approximation algorithms into the construction of a small set that approximates every tuplewise barycenter location. Its main results are as follows.
	\begin{enumerate}
		\item[(i)] We prove the first PTAS for the discrete type-$2$ Wasserstein barycenter problem. For every $\alpha\in(0,1]$, the proposed family of algorithms attains a $(1+\alpha)$ multiplicative guarantee in time polynomial in $(nk)^{1/\alpha}$ and $d$.
		
		\item[(ii)] We provide both randomized sampling and deterministic enumeration variants. Under equal weights, sampling without replacement yields the sharper ratio $1+(k-t)/(t(k-1))$, compared with $1+1/t$ for general weights.
		
		\item[(iii)] We evaluate the algorithms on synthetic and real datasets that include many input measures and large marginal supports. The results show that the reduced-support models are computationally practical and produce high-quality barycenters in regimes where the full MOT formulation is intractable.
		
		\item[(iv)] We extend the support-reduction principle in two directions. For sparse type-$2$ barycenters, the same candidate set yields a $(1+1/t)$ value guarantee when the restricted mixed-integer model is solved to optimality. For type-$1$ barycenters, we derive exact and approximate candidate-support constructions based on coordinatewise weighted medians.
	\end{enumerate}

	The remainder of the paper is organized as follows. Section~\ref{section:reduction} develops the candidate-support reduction framework. Section~\ref{section:W2_withrep} presents randomized and deterministic algorithms for general weights, and Section~\ref{section:W2_improved} sharpens the guarantee under equal weights. Section~\ref{section:experiments} reports the numerical results. Section~\ref{section:extension} studies sparse type-$2$ and type-$1$ extensions, and Section~\ref{section:conclusion} concludes.

	\textbf{Notation.} We write $[k]:=\{1,\ldots,k\}$ and use $\|\cdot\|_2$ for the Euclidean norm. The symbols $\delta_{(\cdot)}$ and $\mathbb{I}_{(\cdot)}$ denote a Dirac measure and an indicator, respectively. We use $(\mathbb{R}_+^n)^{\otimes k}$ for the $k$-fold tensor product of $\mathbb{R}_+^n$, and $\mathbb{P}(S)$ for a probability measure supported on $S$. Finally, $O(\cdot)$ denotes standard asymptotic complexity, whereas $\tilde O(\cdot)$ suppresses logarithmic and polylogarithmic factors in the input size and accuracy.

	\section{Preliminary results and approximation scheme}\label{section:reduction}
	
	This section develops the candidate-support reduction framework that underlies our algorithms and generalizes the construction of \cite{borgwardt2022lp}. We first recall the equivalence between the discrete barycenter problem and multi-marginal optimal transport (MOT), which yields an exact but exponentially large candidate support $S^*$. We then restrict the tuplewise barycenter locations to a smaller set $S$. The resulting dual problem admits a separation oracle with running time polynomial in $n$, $k$, $d$, and $|S|$, and an optimal primal solution induces an explicit barycenter supported on $S$. Finally, we show that a tuplewise approximation guarantee for $S$ implies the same multiplicative guarantee for the barycenter objective.

	\subsection{Optimal support}
	
	A key ingredient is the standard equivalence between the Wasserstein barycenter problem and an MOT problem; see Proposition~2.1 of \cite{altschuler2022wasserstein} and Section~6 of \cite{matching2010}. Specifically, \eqref{eq:def_barycenter} is equivalent to
	\begin{align}
		v^*=&\min_{\bm{\Pi} \in (\mathbb{R}^n_{+})^{\otimes k}} \langle \bm{C},\bm{\Pi} \rangle,\notag\\
		\text{s.t.}&\sum_{j_1\in [|\Xi_1|]}\ldots \sum_{j_{i-1}\in [|\Xi_{i-1}|]}\sum_{j_{i+1}\in [|\Xi_{i+1}|]}\ldots \sum_{j_k\in [|\Xi_k|]}\Pi_{j_1,\ldots,j_k}=\hat{p}_{ij_i},\forall i \in [k], j_i\in [|\Xi_i|],\label{eq_mot}
	\end{align}
	where $\bm{\Pi}$ is a $k$-way coupling tensor and the transportation cost is given by
	$$C_{\vec{\bm{j}}}=\min_{\bm{w} \in \mathbb{R}^d}\sum_{i\in [k]} \lambda_i\|\hat{\bm{x}}_{i,j_i}-\bm{w}\|_2^2,\qquad \vec{\bm{j}}=(j_1,\ldots,j_{i-1},j_i,j_{i+1},\ldots,j_k).$$ 
	Each entry $\Pi_{j_1,\ldots,j_k}$ is a joint mass assigned to the tuple $(\hat{\bm{x}}_{1,j_1},\ldots,\hat{\bm{x}}_{k,j_k})$, and $C_{\vec{\bm{j}}}$ is the minimum weighted squared distance from that tuple to a common barycenter location $\bm w$. The MOT representation immediately yields a finite exact candidate support consisting of tuplewise minimizers.
	\begin{lemma}\label{lem:optimal_support}
		The candidate support $S^*$ of the optimal barycenter of problem \eqref{eq:def_barycenter} is given by
		\[S^*=\bigcup_{\bm{x}_i \in \Xi_i, i\in [k]}\left\{\arg\min_{\bm{w} \in \mathbb{R}^d}\sum_{i\in [k]} \lambda_i\|\bm{x}_{i}-\bm{w}\|_2^2\right\}.\]
	\end{lemma}
	This characterization appears in the remark following Proposition~2.1 of \cite{altschuler2022wasserstein} and in Section~6 of \cite{matching2010}; it is also consistent with the finiteness results of \cite{anderes2016discrete}. Thus, an optimal barycenter can be chosen with finite support, and every support point is a weighted least-squares minimizer associated with a tuple of input atoms.
	
	By Lemma~\ref{lem:optimal_support}, we have $\lvert S^* \rvert = O(n^k)$, which is exponential in $k$. Using this notation, we may equivalently rewrite the transportation cost in the MOT formulation as
	$$C_{\vec{\bm{j}}}=\min_{\bm{w}\in S^*}\sum_{i\in [k]} \lambda_i\|\hat{\bm{x}}_{i,j_i}-\bm{w}\|_2^2,$$ 
	for each index tuple $\vec{\bm{j}} = (j_1,\ldots,j_k)$. Although an optimal barycenter is finitely supported, the exact candidate set $S^*$ may contain exponentially many points. Our approximation framework therefore replaces $S^*$ with a substantially smaller candidate set.

	\subsection{Dual representation of the MOT formulation and candidate-support reduction}
	Let $\gamma_{ij}$ be the dual variable associated with the marginal constraint for atom $j$ of measure $i$. By LP strong duality, the barycenter value admits the representation:
	\begin{align}
		v^*=\max_{\bm{\gamma}}\left\{\sum_{i\in [k]}\sum_{j\in [|\Xi_i|]} \hat{p}_{ij}\gamma_{ij}:C_{\vec{\bm{j}}} - \sum_{i\in [k]} \gamma_{ij_i}\geq 0, \forall \vec{\bm{j}}\in \otimes_{i\in [k]}[|\Xi_i|]\right\}. \label{eq_MOT-D}
	\end{align}
	The dual has at most $nk$ variables, but it contains one constraint for every $k$-tuple $\vec{\bm{j}}$ and therefore exponentially many constraints. Moreover, evaluating a constraint through the exact representation may require minimizing over the exponentially large set $S^*$. We address both difficulties by replacing $S^*$ with a manageable candidate set $S\subseteq\mathbb{R}^d$. Our goal is to ensure that the transportation cost induced by the reduced support set $S$, for every tuple of input support points, is within a prescribed multiplicative factor of the unrestricted cost, either deterministically or in expectation. 

	For a finite candidate support $S=\{\bm w_\ell\}_{\ell\in[|S|]}\subseteq\mathbb{R}^d$, define the restricted tuple cost
	$$C_{\vec{\bm{j}}}(S)=\min_{\bm{w}\in S}\sum_{i\in [k]} \lambda_i\|\hat{\bm{x}}_{i,j_i}-\bm{w}\|_2^2,$$
	for each multi-index $\vec{\bm{j}} = (j_1,\ldots,j_k)$.
	The corresponding restricted dual value is
	\begin{align}
		v(S)=\max_{\bm{\gamma}}\left\{\sum_{i\in [k]}\sum_{j\in [|\Xi_i|]} \hat{p}_{ij}\gamma_{ij}:C_{\vec{\bm{j}}}(S) - \sum_{i\in [k]} \gamma_{ij_i} \geq 0, \forall \vec{\bm{j}}\in \otimes_{i\in [k]}[|\Xi_i|]\right\}. \label{eq_MOT-D_S_p}
	\end{align}
	Equivalently, $v(S)$ is the optimal barycenter value when all tuplewise barycenter locations are restricted to $S$.

	When $|S|$ is polynomial in the input size, the restricted dual can be solved in polynomial time since its separation problem decomposes across marginals.
	\begin{lemma}\label{lem_separation}
		Given a candidate dual vector $\bm{\gamma}$ for \eqref{eq_MOT-D_S_p}, feasibility of $\bm{\gamma}$ can be checked in time $O(k n |S| d)$.
		Therefore, the restricted MOT dual \eqref{eq_MOT-D_S_p} can be solved in time $\tilde{O}(n^3k^3d|S|)$.
	\end{lemma}
	\begin{proof} We observe that the fact that $\bm{\gamma}$ is feasible to \eqref{eq_MOT-D_S_p} is equivalent to
		\[\min_{\vec{\bm{j}}\in \otimes_{i\in [k]}[|\Xi_i|]}C_{\vec{\bm{j}}}(S) - \sum_{i\in [k]} \gamma_{ij_i} \geq 0,\]
		which is equivalent to the following two-stage minimization problem:
		\[\min_{\ell\in [|S|]}\min_{\vec{\bm{j}}\in \otimes_{i\in [k]}[|\Xi_i|]}\sum_{i\in [k]} \lambda_i\|\bm{w}_{\ell}-\hat{\bm{x}}_{ij_i}\|_2^2 - \sum_{i\in [k]} \gamma_{ij_i} \geq 0.\]
		Note that the inner minimization problem can be decomposed for each $i\in [k]$ when $\ell\in [|S|]$ is fixed. Let us define the cost $c_{ij_i\ell}=\lambda_i\|\bm{w}_{\ell}-\hat{\bm{x}}_{ij_i}\|_2^2-\gamma_{ij_i}$. Then for each $i \in [k]$, let  $j_i^*(\ell)\in \arg\min_{j_i}c_{ij_i\ell}$. Next, we solve the outer minimization by picking the $\ell^*$ such that
		$\ell^*\in \arg\min_{\ell\in [|S|]}\sum_{i\in [k]}c_{ij_i^*(\ell)\ell}.$
		
		Computing each cost coefficient $c_{ij_i\ell}$ requires $O(d)$ time. Therefore, computing all such coefficients costs $O(kn|S|d)$ in total, solving the inner minimization requires $O(kn|S|)$, and solving the outer minimization takes $O(k|S|)$. In total, the complexity of the separation oracle is $O(kn|S|d)$. By the complexity result of the ellipsoid method with a separation oracle (Section 3.1 in \cite{grotschel2012geometric}), the overall complexity for solving the restricted MOT dual \eqref{eq_MOT-D_S_p} is $\tilde{O}((nk)^2kn|S|d)=\tilde{O}(k^3n^3|S|d)$.
	\end{proof}
	
	An optimal restricted MOT solution also yields an explicit barycenter supported on $S$.
	\begin{lemma}\label{lem_distr}
		Fix a candidate support $S \subseteq \mathbb{R}^d$. Let $\bm{\Pi}^*(S)$ be an optimal primal solution of the MOT formulation \eqref{eq_mot} in which the cost $C_{\vec{\bm{j}}}$ is restricted to $C_{\vec{\bm{j}}}(S)$. For each $\vec{\bm{j}}\in \otimes_{i\in [k]}[|\Xi_i|]$ such that $\Pi_{\vec{\bm{j}}}^*(S) > 0$, choose
		$$\bm{w}_{\vec{\bm{j}}}\in \arg\min_{\bm{w}\in S}\sum_{i\in [k]} \lambda_i\|\hat{\bm{x}}_{i,j_i}-\bm{w}\|_2^2,\qquad m_{\bm{w}}=\sum_{\vec{\bm{j}}\in\otimes_{i\in [k]}[|\Xi_i|]}\Pi_{
			\vec{\bm{j}}
		}^*(S)\mathbb{I}_{\{\bm{w}_{\vec{\bm{j}}}=\bm{w}\}},\forall \bm{w}\in S, $$
		where $\mathbb{I}_{(\cdot)}$ is the indicator function. 
		Then the measure $\mathbb{P}(S)=\sum_{\bm{w}\in S}m_{\bm{w}}\delta_{\tilde{\bm{w}}=\bm{w}}$ 
		is a feasible barycenter supported on $S$, where $\delta(\cdot)$ is the Dirac measure and the random variable $\tilde{\bm{w}}$ follows the distribution of $\mathbb{P}(S)$.
	\end{lemma}
	The tensor $\bm\Pi^*(S)$ specifies the joint mass assigned to each tuple of input atoms. Assigning that mass to a best-fitting point $\bm w_{\vec{\bm j}}\in S$ and aggregating over tuples produces the measure $\mathbb P(S)$. We next quantify how the restriction from $S^*$ to a smaller candidate set $\tilde S$ affects the optimal value.
	\begin{thm}\label{thm:ratio_support_reduction}
		Let $\tilde S \subseteq \mathbb{R}^d$ be a (possibly random) candidate support set. Suppose there exists $\alpha \ge 0$ such that, for any choice of points $\bm{x}_i \in \Xi_i$ for all $i \in [k]$, we have
		\begin{align}
			\mathbb{E}_{\tilde{S}}\left[\min_{\bm{w} \in \tilde{S}}\sum_{i\in [k]}\lambda_i\|\bm{w}-\bm{x}_i\|_2^2\right]\leq (1+\alpha)\min_{\bm{w} \in \mathbb{R}^d}\sum_{i\in [k]}\lambda_i\|\bm{w}-\bm{x}_i\|_2^2. \label{eq:ratio_support_reduction}
		\end{align}
		Then the expected objective value of the restricted barycenter problem satisfies $\mathbb{E}_{\tilde{S}}[v(\tilde{S})] \leq (1+\alpha) v(S^*)$.
	\end{thm}
	\begin{proof}
		We first observe that in the MOT formulation \eqref{eq_mot}, replacing transportation cost $\bm{C}$ by the approximate one $\bm{C}(\tilde{S})$, the value of $v(\tilde{S})$ can be obtained by solving the following restricted MOT problem:
		\begin{align}
			&v(\tilde{S})=\min_{\bm{\Pi} \in (\mathbb{R}^n_{+})^{\otimes k}} \langle \bm{C}(\tilde{S}),\bm{\Pi} \rangle,\notag\\
			\text{s.t.}&\sum_{j_1\in [|\Xi_1|]}\ldots \sum_{j_{i-1}\in [|\Xi_{i-1}|]}\sum_{j_{i+1}\in [|\Xi_{i+1}|]}\ldots \sum_{j_k\in [|\Xi_k|]}\Pi_{j_1,\ldots,j_{i-1},j_i,j_{i+1},\ldots,j_k}=\hat{p}_{ij_i},\forall i \in [k], j_i\in [|\Xi_i|].\label{eq_mot_tilde_S}
		\end{align}
		Suppose that $\bm{\Pi}^*$ is an optimal solution of the MOT formulation \eqref{eq_mot}, which is also feasible to the restricted MOT formulation \eqref{eq_mot_tilde_S}. Thus,
		\begin{align*}
			v(\tilde{S})\leq \langle \bm{C}(\tilde{S}),\bm{\Pi}^* \rangle=\sum_{\vec{\bm{j}}\in \otimes_{i\in [k]}[|\Xi_i|]}C_{\vec{\bm{j}}}(\tilde{S})\Pi^*_{\vec{\bm{j}}},
		\end{align*}
		where $\vec{\bm{j}}=(j_1,\ldots,j_k)$. Taking expectation on both sides of the inequality, we have
		\begin{align*}
			\mathbb{E}_{\tilde{S}}\left[v(\tilde{S})\right]\leq \sum_{\vec{\bm{j}}\in \otimes_{i\in [k]}[|\Xi_i|]}\mathbb{E}_{\tilde{S}}\left[C_{\vec{\bm{j}}}(\tilde{S})\Pi^*_{\vec{\bm{j}}}\right] \leq (1+\alpha)\sum_{\vec{\bm{j}}\in \otimes_{i\in [k]}[|\Xi_i|]}\mathbb{E}_{\tilde{S}}\left[C_{\vec{\bm{j}}}(S^*)\Pi^*_{\vec{\bm{j}}}\right]:=(1+\alpha)v(S^*),
		\end{align*}
		where the second inequality is due to the condition \eqref{eq:ratio_support_reduction} and the fact that $\Pi^*_{\vec{\bm{j}}}\geq 0$ for any $\vec{\bm{j}}$.
	\end{proof}
	\begin{remark}
		If the set $\tilde S$ is deterministic, then the same guarantee holds directly under condition~\eqref{eq:ratio_support_reduction} (i.e., the expectation over $\tilde S$ is not needed).
	\end{remark}
	Theorem~\ref{thm:ratio_support_reduction} reduces the global barycenter approximation problem to a local one: it suffices to approximate the optimal one-center cost for every tuple of input atoms. The next two sections construct candidate supports satisfying \eqref{eq:ratio_support_reduction}, either deterministically or in expectation.
	
	\section{Approximation algorithms for the Wasserstein barycenter problem with general weights} \label{section:W2_withrep}
	Under squared Euclidean loss, the tuplewise minimizer in \Cref{lem:optimal_support} is the weighted centroid. Hence, an optimal barycenter can be supported on
	\begin{align*}
		S^*:=\left\{\sum_{i\in [k]}\lambda_i\bm{x}_i: \bm{x}_i \in \Xi_i, \forall i\in [k]\right\}.
	\end{align*}
	This set contains at most $\prod_{i\in[k]}|\Xi_i|=O(n^k)$ points and is therefore exponential in $k$. For general weights $\bm\lambda\in\Delta_k$, we construct smaller supports by averaging $t$ sampled input atoms. Sampling with replacement preserves the weighted centroid in expectation, which leads to both a randomized algorithm and a deterministic enumeration counterpart with the same approximation ratio.
	
	\subsection{Sampling with replacement}
	Draw $T_1,\ldots,T_t\in[k]$ independently according to $\bm\lambda$, so that
	\begin{align*}
		\textrm{Prob}\left[T_j=i\right]=\lambda_i,  \forall  i\in[k],\ j\in[t].
	\end{align*}
	Let $\mathcal T$ be the resulting multiset. We construct
	\begin{align*}
		S_1^{\mathcal{T}}=\left\{\frac{1}{t}\sum_{i \in [t]}\bm{x}_{T_i}:\bm{x}_{T_i} \in \Xi_{T_i},\forall i \in [t]\right\}. 
	\end{align*}
	Thus, every candidate atom is the average of $t$ input atoms, with one atom chosen from each sampled marginal. The support size is at most $n^t$, and hence polynomial in $n$ for fixed $t$. Solving the restricted MOT problem over $S_1^{\mathcal T}$ and applying \Cref{lem_distr} produces the approximate barycenter $\mathbb P(S_1^{\mathcal T})$.
	\begin{algorithm}[htbp]
		\caption{Sampling-with-replacement algorithm for the Wasserstein barycenter problem}
		\label{alg:subset_sampling_with_rep}
		\begin{algorithmic}[1]
			\State \textbf{Input}: Probability measures $\{\mathbb{P}_i\}_{i\in[k]}$ with supports $\{\Xi_i\}_{i\in[k]}$ and an integer $t\in[k]$.
			\State Independently draw $T_1,\ldots,T_t\in[k]$ with $\textrm{Prob}[T_j=i]=\lambda_i$, and let $\mathcal T$ be the resulting multiset.
			\State $S_1^{\mathcal{T}}:=\left\{\frac{1}{t}\sum_{i \in [t]}\bm{x}_{T_i}: \bm{x}_{T_i} \in \Xi_{T_i} \right\}$.
			\State Solve the restricted MOT dual \eqref{eq_MOT-D_S_p} and recover a corresponding optimal primal solution $\bm{\Pi}^*(S_1^{\mathcal{T}})$.
			\State \textbf{Output}: Construct and return the approximate barycenter $\mathbb{P}(S_1^{\mathcal{T}})$ using \Cref{lem_distr}.
		\end{algorithmic}
	\end{algorithm}
	
	The approximation guarantee is formally stated in the following theorem.
	\begin{thm}\label{thm:W2_sampling_with_rep}
		Let $t \in [k]$ be a fixed constant, and let $S_1^{\mathcal{T}}$ be the candidate support returned by \Cref{alg:subset_sampling_with_rep}. Then
		\begin{align*}
			\mathbb{E}_{\mathcal{T}}\left[v\left(S_1^{\mathcal{T}}\right)\right] \leq \left(1+\frac{1}{t}\right)v^*.
		\end{align*}
	\end{thm}
	\begin{proof}
		\begin{subequations}
			By Theorem \ref{thm:ratio_support_reduction}, it is sufficient to show that for any fixed realization $\bm{x}_i \in \Xi_i$ for all $i\in [k]$ and $\bm{c}:=\sum_{i\in [k]} \lambda_i\bm{x}_i$, we have
			\begin{align}
				\mathbb{E}_{\mathcal{T}}\left[\min_{\bm{s} \in S_1^{\mathcal{T}}}\sum_{i\in [k]}\lambda_i\|\bm{s}-\bm{x}_i\|_2^2\right] \leq \left(1+\frac{1}{t}\right)\sum_{i\in [k]}\lambda_i\|\bm{c}-\bm{x}_i\|_2^2.\label{eq_exp_condition}
			\end{align}
			Since we consider a particular realization from each marginal distribution, for each random sample $T_j$, we also consider the same realization; i.e., we let $\bm{x}_{T_j}=\sum_{i\in [k]}\bm{x}_i \mathbb{I}_{\{T_j=i\}}$.
			
			Note that $\bm{c} \in \arg\min_{\bm{s} \in \mathbb{R}^d}\sum_{i \in [k]}\lambda_i\|\bm{s}-\bm{x}_i\|_2^2$. Therefore, to approximate the solution of $\min_{\bm{s} \in S_1^{\mathcal{T}}}\sum_{i\in [k]}\lambda_i\|\bm{s}-\bm{x}_i\|_2^2$, we consider $\bm{s}_c \in \arg\min_{\bm{s} \in S_1^{\mathcal{T}}}\|\bm{s}-\bm{c}\|_2^2$; i.e., the point in $S_1^{\mathcal{T}}$ that is closest to $\bm{c}$. We have
			\begin{align}
				&\sum_{i\in [k]}\lambda_i\|\bm{s}_c-\bm{x}_i\|_2^2=\bm{s}_c^{\top}\bm{s}_c-2\bm{c}^{\top}\bm{s}_c+\sum_{i\in [k]}\lambda_i \bm{x}_i^{\top}\bm{x}_i
				=\sum_{i\in [k]}\lambda_i(\|\bm{s}_c-\bm{c}\|_2^2+\|\bm{c}-\bm{x}_i\|_2^2),\label{eq:l2_expansion}
			\end{align}
			where the equalities follow by direct expansion and $\sum_{i \in [k]}\lambda_i=1$.
			
			Since $\bm{s}_c\in S_1^{\mathcal{T}}$, by \eqref{eq:l2_expansion}, the left-hand side of \eqref{eq_exp_condition} can be upper bounded as
			\begin{align}
				\mathbb{E}_{\mathcal{T}}\left[\min_{\bm{s} \in S_1^{\mathcal{T}}}\sum_{i\in [k]}\lambda_i\|\bm{s}-\bm{x}_i\|_2^2\right] \leq \sum_{i\in [k]}\lambda_i\|\bm{c}-\bm{x}_i\|_2^2+ \mathbb{E}_{\mathcal{T}}\left[\|\bm{s}_c-\bm{c}\|_2^2\right].
				\label{eq_exp_condition_v2}
			\end{align}
			It remains to bound $\mathbb{E}_{\mathcal{T}}\left[\|\bm{s}_c-\bm{c}\|_2^2\right]$. Since $\frac{\sum_{i\in \mathcal{T}}\bm{x}_i}{t}\in S_1^{\mathcal{T}}$, we have
			\begin{align}\label{eq:min vs expectation_sample_with_rep}
				\mathbb{E}_{\mathcal{T}}\left[\|\bm{s}_c-\bm{c}\|_2^2\right]=\mathbb{E}_{\mathcal{T}}\left[\min_{\bm{s} \in S_1^{\mathcal{T}}}\|\bm{s}-\bm{c}\|_2^2\right] \leq \mathbb{E}_{\mathcal{T}}\left[\left\|\frac{\sum_{i\in \mathcal{T}}\bm{x}_i}{t}-\bm{c}\right\|_2^2\right].
			\end{align}
			Since $\mathcal{T}=\{T_i\}_{i\in [k]}$ and $T_i$'s are i.i.d., we have
			\begin{align*}
				\mathbb{E}\left[\frac{\sum_{i \in \mathcal{T}}\bm{x}_i}{t}\right]=\frac{\sum_{j \in [t]}\sum_{i\in [k]}\mathrm{{Prob}}\left[T_j=i\right]\bm{x}_{i}}{t}=\frac{t\sum_{i\in [k]}\lambda_i\bm{x}_i}{t}=\sum_{i\in [k]}\lambda_i\bm{x}_i=\bm{c}.
			\end{align*}
			Therefore, we have
			\begin{align}\label{eq:min vs expectation_sample_with_rep_v2}
				\mathbb{E}_{\mathcal{T}}\left[\|\bm{s}_c-\bm{c}\|_2^2\right]\leq \mathbb{E}_{\mathcal{T}}\left[\left\|\frac{\sum_{i \in \mathcal{T}}\bm{x}_{i}}{t}-\bm{c}\right\|_2^2\right]=\frac{1}{t}\mathbb{E}_{T_1}\left[\left\|\bm{x}_{T_1}-\bm{c}\right\|_2^2\right]=\frac{1}{t}\sum_{i\in [k]}\lambda_i\|\bm{x}_i-\bm{c}\|_2^2.
			\end{align}
			
			Plugging the upper bound of \eqref{eq:min vs expectation_sample_with_rep_v2} into \eqref{eq_exp_condition_v2}, we obtain that
			\begin{align*}
				\mathbb{E}_{\mathcal{T}}\left[\min_{\bm{s} \in S_1^{\mathcal{T}}}\sum_{i\in [k]}\lambda_i\|\bm{s}-\bm{x}_i\|_2^2\right] \leq \left(1+\frac{1}{t}\right)\sum_{i\in [k]}\lambda_i\|\bm{c}-\bm{x}_i\|_2^2.
			\end{align*}
			This completes the proof.
		\end{subequations}
	\end{proof}
	\begin{remark}
		Choosing $t=\lceil1/\alpha\rceil$ gives an expected $(1+\alpha)$-approximation. Since $|S_1^{\mathcal T}|=O(n^t)$, \Cref{lem_separation} yields the running-time bound $\tilde O(n^{3+\lceil1/\alpha\rceil}k^3d)$. Thus, for every fixed $\alpha>0$, the randomized algorithm runs in polynomial time.
	\end{remark}
	
	\subsection{Deterministic counterpart based on multiset enumeration}
	The randomized construction can be derandomized by enumerating all multisets of $t$ marginal indices and collecting every corresponding average. Define
	\begin{align}
		S_1^{t}:=\left\{\frac{1}{t}\sum_{i \in [t]}\bm{x}_{T_i}: \bm{x}_{T_i} \in \Xi_{T_i}, T_i \in [k],i\in[t]\right\}. \label{eq:candidate_support_with_rep_deterministic}
	\end{align}
	Solving the restricted MOT problem over $S_1^t$ and applying \Cref{lem_distr} yields the deterministic approximate barycenter $\mathbb P(S_1^t)$.
	
	\begin{algorithm}[htbp]
		\caption{Multiset-enumeration algorithm for the Wasserstein barycenter problem}
		\label{alg:subset_enumeration_with_rep}
		\begin{algorithmic}[1]
			\State \textbf{Input}: Probability measures $\{\mathbb{P}_i\}_{i\in[k]}$ with supports $\{\Xi_i\}_{i\in[k]}$ and an integer $t\in[k]$.
			\State $S_1^{t}:=\left\{\frac{1}{t}\sum_{i \in [t]}\bm{x}_{T_i}: \bm{x}_{T_i} \in \Xi_{T_i}, T_i \in [k],i\in[t]\right\}$.
			\State Solve the restricted MOT dual \eqref{eq_MOT-D_S_p} and recover a corresponding optimal primal solution $\bm{\Pi}^*(S_1^{t})$.
			\State \textbf{Output}: Construct and return the approximate barycenter $\mathbb{P}(S_1^{t})$ using \Cref{lem_distr}.
		\end{algorithmic}
	\end{algorithm}
	
	Since $S_1^t$ contains the candidate support generated by every realization of $\mathcal T$, for any reference point $\bm c$, we have
	\begin{align*}
		\min_{\bm{s} \in S_1^{t}} \|\bm{s}-\bm{c}\|_2^2 \leq \mathbb{E}_{\mathcal{T}}\left[\min_{\bm{s} \in S_1^{\mathcal{T}}}\|\bm{s}-\bm{c}\|_2^2\right] \leq \mathbb{E}_{T}\left[\left\|\frac{\sum_{i\in \mathcal{T}}\bm{x}_i}{t}-\bm{c}\right\|_2^2\right],
	\end{align*}
	The same argument as in Theorem~\ref{thm:W2_sampling_with_rep} therefore gives the deterministic guarantee below.
	\begin{thm}\label{thm:W2_enumeration_with_rep}
		Let $t \in [k]$ be a fixed constant, and let $S_1^{t}$ be the candidate support returned by \Cref{alg:subset_enumeration_with_rep}. Then
		\begin{align*}
			v(S_1^{t}) \leq \left(1+\frac{1}{t}\right)v^*,
		\end{align*}
		and the support size satisfies $|S_1^{t}|=O((nk)^t)$.
	\end{thm}
	\begin{proof}
		According to the proof of Theorem \ref{thm:W2_sampling_with_rep}, it is sufficient to show that
		\begin{align*}
			\min_{\bm{s} \in S_1^{t}}\sum_{i\in [k]}\lambda_i\|\bm{s}-\bm{x}_i\|_2^2 \leq \mathbb{E}_{\mathcal{T}}\left[\min_{\bm{s} \in S_1^{\mathcal{T}}}\sum_{i\in [k]}\lambda_i\|\bm{s}-\bm{x}_i\|_2^2\right],
		\end{align*}
		holds for any $\bm{x}_i \in \Xi_i$ and $i\in [k]$. This is indeed true since by definition, we have $S_1^{\hat{T}}\subseteq S_1^{t}$ for any realization $\hat{T}$ of random set $\mathcal{T}$. The approximation ratio then follows from Theorem \ref{thm:W2_sampling_with_rep}.
		
		Finally, note that there are at most $\binom{k+t-1}{t} n^t$ points in the set $S_1^{t}$ by stars and bars method. Therefore, the size of $S_1^{t}$ is $O(nk)^{t}$.
	\end{proof}
	\begin{remark}
		For $t=1$, $S_1^t$ is the union of the input supports, so Algorithm~\ref{alg:subset_enumeration_with_rep} reduces to the $2$-approximation of \cite{borgwardt2022lp}. More generally, choosing $t=\lceil1/\alpha\rceil$ gives a deterministic $(1+\alpha)$-approximation. Since $|S_1^t|=O((nk)^t)$, the resulting running time is $\tilde O(n^{3+\lceil1/\alpha\rceil}k^{3+\lceil1/\alpha\rceil}d)$.
	\end{remark}

	\section{Approximation algorithms for the Wasserstein barycenter problem with equal weights}\label{section:W2_improved}
	We now assume equal weights, $\lambda_i=1/k$ for all $i\in[k]$. Sampling marginals without replacement reduces the variance of the sampled centroid and consequently improves the approximation ratio relative to the general-weight construction.

	\subsection{Sampling without replacement}
	Uniformly sample a subset $\mathcal T\subseteq[k]$ of cardinality $t$ and define
	\begin{align*}
		S_2^{\mathcal{T}}:=\left\{\frac{1}{t}\sum_{i\in \mathcal{T}}\bm{x}_i: \bm{x}_i \in \Xi_i \right\}.
	\end{align*}
	The restricted MOT problem over $S_2^{\mathcal T}$ yields an approximate barycenter through \Cref{lem_distr}.
	\begin{algorithm}[htbp]
		\caption{Sampling-without-replacement algorithm for the Wasserstein barycenter problem}
		\label{alg:subset_sampling}
		\begin{algorithmic}[1]
			\State \textbf{Input}: Probability measures $\{\mathbb{P}_i\}_{i\in[k]}$ with supports $\{\Xi_i\}_{i\in[k]}$ and an integer $t\in[k]$.
			\State Uniformly sample a subset $\mathcal T\subseteq[k]$ with $|\mathcal T|=t$.
			\State $S_2^{\mathcal{T}}:=\left\{\frac{1}{t}\sum_{i\in \mathcal{T}}\bm{x}_i: \bm{x}_i \in \Xi_i\right\}$.
			\State Solve the restricted MOT dual \eqref{eq_MOT-D_S_p} and recover a corresponding optimal primal solution $\bm{\Pi}^*(S_2^{\mathcal{T}})$.
			\State \textbf{Output}: Construct and return the approximate barycenter $\mathbb{P}(S_2^{\mathcal{T}})$ using \Cref{lem_distr}.
		\end{algorithmic}
	\end{algorithm}
	
	The approximation guarantee is formally stated in the following theorem.
	\begin{thm}\label{thm:W2_sampling}
		Let $t \in [k]$ be a fixed constant, and let $S_2^{\mathcal{T}}$ be the candidate support returned by \Cref{alg:subset_sampling}. Then
		\begin{align*}
			\mathbb{E}_{\mathcal{T}}\left[v(S_2^{\mathcal{T}})\right] \leq \left(1+\frac{k-t}{t(k-1)}\right)v^*.
		\end{align*}
	\end{thm}
	\begin{subequations}
		\begin{proof}
			By Theorem \ref{thm:ratio_support_reduction}, we only need to show that for any fixed realization $\bm{x}_i \in \Xi_i$, $i\in [k]$ and $\bm{c}:=\frac{1}{k}\sum_{i\in [k]}\bm{x}_i$, we must have
			\begin{align}\label{eq_equivalent_w2_equal}
				\mathbb{E}_{\mathcal{T}}\left[\min_{\bm{s} \in S_2^{\mathcal{T}}}\frac{1}{k}\sum_{i\in [k]}\|\bm{s}-\bm{x}_i\|_2^2\right] \leq \left(1+\frac{k-t}{t(k-1)}\right)\frac{1}{k}\sum_{i\in [k]}\|\bm{c}-\bm{x}_i\|_2^2.
			\end{align}
			Note that $\bm{c} \in \arg\min_{\bm{s} \in \mathbb{R}^d}\frac{1}{k}\sum_{i \in [k]}\|\bm{s}-\bm{x}_i\|_2^2$. Therefore, to approximate the solution of $\min_{\bm{s} \in S_2^{\mathcal{T}}}\frac{1}{k}\sum_{i\in [k]}\|\bm{s}-\bm{x}_i\|_2^2$, we consider $\bm{s}_c \in \arg\min_{\bm{s} \in S_2^{\mathcal{T}}}\|\bm{s}-\bm{c}\|_2^2$, which is the closest point in $S_2^{\mathcal{T}}$ to $\bm{c}$. Similar to the proof of Theorem \ref{thm:W2_sampling_with_rep}, we have
			\begin{align*}
				\frac{1}{k}\sum_{i\in [k]}\|\bm{s}_{\bm{c}}-\bm{x}_i\|_2^2
				=\frac{1}{k}\sum_{i\in [k]}(\|\bm{s}_{\bm{c}}-\bm{c}\|_2^2+\|\bm{c}-\bm{x}_i\|_2^2).
			\end{align*}
			Thus, the left-hand side can be upper bounded by
			\begin{align}\label{eq_equivalent_w2_equal2}
				\mathbb{E}_{\mathcal{T}}\left[\min_{\bm{s} \in S_2^{\mathcal{T}}}\frac{1}{k}\sum_{i\in [k]}\|\bm{s}-\bm{x}_i\|_2^2\right] \leq 
				\frac{1}{k}\sum_{i\in [k]}\|\bm{c}-\bm{x}_i\|_2^2+\mathbb{E}_{\mathcal{T}}\left[\min_{\bm{s} \in S_2^{\mathcal{T}}}\|\bm{s}-\bm{c}\|_2^2\right].
			\end{align}
			
			It remains to bound $\mathbb{E}_{\mathcal{T}}\left[\min_{\bm{s} \in S_2^{\mathcal{T}}}\|\bm{s}-\bm{c}\|_2^2\right]$. 
			For notational convenience, we let $\hat{\bm{\mu}}_{\mathcal{T}}=\frac{\sum_{i\in \mathcal{T}}\bm{x}_i}{t}$. Since $\hat{\bm{\mu}}_{\mathcal{T}}\in S_2^{\mathcal{T}}$, we have
			\begin{align}\label{eq:min_vs_expectation}
				\mathbb{E}_{\mathcal{T}}\left[\min_{\bm{s} \in S_2^{\mathcal{T}}}\|\bm{s}-\bm{c}\|_2^2\right] \leq \mathbb{E}_{\mathcal{T}}\left[\left\|\hat{\bm{\mu}}_{\mathcal{T}}-\bm{c}\right\|_2^2\right].
			\end{align}
			
			Note that
			\begin{align*}
				\mathbb{E}_{\mathcal{T}}\left[\hat{\bm{\mu}}_{\mathcal{T}}\right] &=\sum_{\hat{T}\in {[k]}\choose t}\mathrm{Prob}[\mathcal{T}=\hat{T}]\frac{\sum_{i \in \hat{T}}\bm{x}_i}{t} =\frac{1}{t\binom{k}{t}} \sum_{\hat{T}\in {[k]}\choose t}\sum_{i \in \hat{T}}\bm{x}_i\\
				&=\frac{1}{t\binom{k}{t}} \sum_{i\in [k]}\bm{x}_i\binom{k-1}{t-1}=\frac{1}{k}\sum_{i\in [k]}\bm{x}_i=\bm{c}.
			\end{align*}
			In addition, we also have
			\begin{align*}
				\text{Cov}\left[\hat{\bm{\mu}}_{\mathcal{T}}\right]=\frac{1}{t^2}\sum_{i,j \in \mathcal{T}}(\bm{x}_i-\bm{c})(\bm{x}_j-\bm{c})^{\top}\mathrm{Prob}\left[i,j \in \mathcal{T}\right],
			\end{align*}
			where 
			\begin{align*}
				\mathrm{Prob}[i,j \in \mathcal{T}]=\begin{cases}
					\frac{t}{k}, &\text{if }i=j,\\
					\frac{\binom{k-2}{t-2}}{\binom{k}{t}}=\frac{t(t-1)}{k(k-1)}, &\text{if }i \neq j,
				\end{cases}
			\end{align*}
			by direct computation.
			
			Therefore, we further have
			\begin{align*}
				\text{Cov}\left[\hat{\bm{\mu}}_{\mathcal{T}} \right]=&\frac{1}{t^2}\sum_{i\in [k]}(\bm{x}_i-\bm{c})(\bm{x}_i-\bm{c})^{\top}\frac{t}{k}+\frac{1}{t^2}\sum_{i \neq j \in [k]}(\bm{x}_i-\bm{c})(\bm{x}_j-\bm{c})^{\top}\frac{t(t-1)}{k(k-1)},
			\end{align*}
			which implies that
			\begin{align}
				\text{Var}\left[\hat{\bm{\mu}}_{\mathcal{T}} \right]=&\mathbb{E}_{\mathcal{T}}\left[\|\hat{\bm{\mu}}_{\mathcal{T}}-\bm{c}\|_2^2\right]\notag\\
				=&\frac{1}{tk}\sum_{i\in [k]}\|\bm{c}-\bm{x}_i\|_2^2+\frac{t-1}{tk(k-1)}\sum_{i\neq j\in [k]}(\bm{x}_i-\bm{c})^\top(\bm{x}_j-\bm{c})\notag\\
				=&\frac{1}{tk}\sum_{i\in [k]}\|\bm{c}-\bm{x}_i\|_2^2-\frac{t-1}{tk(k-1)}\sum_{i\in [k]}\|\bm{c}-\bm{x}_i\|_2^2\notag\\
				=&\frac{k-t}{t(k-1)}\frac{1}{k}\sum_{i\in [k]}\|\bm{c}-\bm{x}_i\|_2^2\label{eq_equivalent_w2_equal3}
			\end{align}
			where the third equality is due to $\sum_{i\in [k]}(\bm{x}_i-\bm{c})=0$.
			
			Combining \eqref{eq_equivalent_w2_equal2}-\eqref{eq_equivalent_w2_equal3}, we have
			\begin{align*}
				\mathbb{E}_{\mathcal{T}}\left[\min_{\bm{s} \in S_2^{\mathcal{T}}}\frac{1}{k}\sum_{i\in [k]}\|\bm{s}-\bm{x}_i\|_2^2\right] \leq  \frac{k-t}{tk(k-1)}\sum_{j \in [k]}\|\bm{c}-\bm{x}_j\|_2^2+\frac{1}{k}\sum_{i\in [k]}\|\bm{c}-\bm{x}_i\|_2^2=\left(1+\frac{k-t}{t(k-1)}\right)\frac{1}{k}\sum_{i\in [k]}\|\bm{c}-\bm{x}_i\|_2^2.
			\end{align*}
		\end{proof}
	\end{subequations}
	
	\begin{remark}
		A $(1+\alpha)$ guarantee is obtained by choosing
		\[
		t=\left\lceil\frac{k}{1+\alpha(k-1)}\right\rceil,
		\]
		which ensures $(k-t)/(t(k-1))\leq\alpha$. Since $|S_2^{\mathcal T}|=O(n^t)$, the corresponding running time is $\tilde O(n^{3+t}k^3d)=\tilde O(n^{3+\lceil 1/\alpha\rceil}k^3d)$.
	\end{remark}
	
	\subsection{Deterministic counterpart based on subset enumeration}
	Enumerating all size-$t$ subsets of marginals gives the deterministic candidate support
	\begin{align*}
		S_2^{t}:=\left\{\frac{1}{t}\sum_{i \in \mathcal{T}}\bm{x}_i: \bm{x}_i \in \Xi_i, \mathcal{T} \subseteq [k],|\mathcal{T}|=t\right\}.
	\end{align*}
	We solve the restricted MOT problem over $S_2^t$ and recover $\mathbb P(S_2^t)$ through \Cref{lem_distr}.
	
	\begin{algorithm}[htbp]
		\caption{Subset-enumeration algorithm for the Wasserstein barycenter problem}
		\label{alg:subset_enumeration}
		\begin{algorithmic}[1]
			\State \textbf{Input}: Probability measures $\{\mathbb{P}_i\}_{i\in[k]}$ with supports $\{\Xi_i\}_{i\in[k]}$ and an integer $t\in[k]$.
			\State $S_2^{t}:=\left\{\frac{1}{t}\sum_{i \in \mathcal{T}}\bm{x}_i: \bm{x}_i \in \Xi_i, \mathcal{T} \subseteq [k],|\mathcal{T}|=t\right\}$.
			\State Solve the restricted MOT dual \eqref{eq_MOT-D_S_p} and recover a corresponding optimal primal solution $\bm{\Pi}^*(S_2^{t})$.
			\State \textbf{Output}: Construct and return the approximate barycenter $\mathbb{P}(S_2^{t})$ using \Cref{lem_distr}.
		\end{algorithmic}
	\end{algorithm}

	Since $S_2^t$ contains the support generated by every realization of $\mathcal T$, for any reference point $\bm c$ we have
	\begin{align*}
		\min_{\bm{s} \in S_2^{t}} \|\bm{s}-\bm{c}\|_2^2 \leq \mathbb{E}_{\mathcal{T}}\left[\min_{\bm{s} \in S_2^{\mathcal{T}}}\|\bm{s}-\bm{c}\|_2^2\right] \leq \mathbb{E}_{\mathcal{T}}\left[\left\|\frac{\sum_{i \in \mathcal{T}}\bm{x}_i}{t}-\bm{c}\right\|_2^2\right].
	\end{align*}
	Applying the argument of Theorem~\ref{thm:W2_sampling} yields the following deterministic guarantee.
	\begin{thm}\label{thm:W2_enumeration}
		Let $t \in [k]$ be a fixed constant, and let $S_2^{t}$ be the candidate support returned by \Cref{alg:subset_enumeration}. Then
		\begin{align*}
			v(S_2^{t}) \leq \left(1+\frac{k-t}{t(k-1)}\right)v^*,
		\end{align*}
		and the support size satisfies $|S_2^{t}|=O((nk)^t)$
	\end{thm}
	\begin{proof}
		According to the proof of Theorem \ref{thm:W2_sampling}, it is sufficient to show that
		\begin{align*}
			\min_{\bm{s} \in S_2^{t}}\sum_{i\in [k]}\lambda_i\|\bm{s}-\bm{x}_i\|_2^2 \leq \mathbb{E}_{\mathcal{T}}\left[\min_{\bm{s} \in S_2^{\mathcal{T}}}\sum_{i\in [k]}\lambda_i\|\bm{s}-\bm{x}_i\|_2^2\right],
		\end{align*}
		holds for any $\bm{x}_i \in \Xi_i$ and any $i\in [k]$. This is indeed true since by definition, we have $S_2^{\hat{T}}\subseteq S_2^{t}$ for any realization $\hat{T}$ of random set $\mathcal{T}$. The approximation ratio then follows from Theorem \ref{thm:W2_sampling}.
		
		Finally, note that there are $\binom{k}{t}n^t$ points in the set $S_2^{t}$. Therefore, the size of $S_2^{t}$ is $O((nk)^{t})$.
	\end{proof}
	
	\begin{remark}
		Since
		\begin{align*}
			1+\frac{k-t}{t(k-1)} \leq 1+\frac{1}{t},
		\end{align*}
		the equal-weight algorithm is never worse than the general-weight construction for the same $t$ and is strictly better whenever $1<t<k$.
	\end{remark}
	
	\section{Numerical experiments}\label{section:experiments}
	We evaluate the computational efficiency and solution quality of the proposed algorithms. The experiments are designed to (i) assess the accuracy--complexity tradeoff controlled by $t$; (ii) demonstrate scalability on instances for which the extensive MOT formulation is intractable; and (iii) compare against the $2$-approximation of \cite{borgwardt2022lp}, which coincides with our deterministic construction at $t=1$. All methods were implemented in Python with Gurobi~11.0.3 and run on a MacBook Pro equipped with an Apple M4 Pro chip and 48~GB of memory.

	\subsection{Synthetic nested ellipse dataset}
\begin{figure}[htbp]
		\centering
		\includegraphics[width=\linewidth]{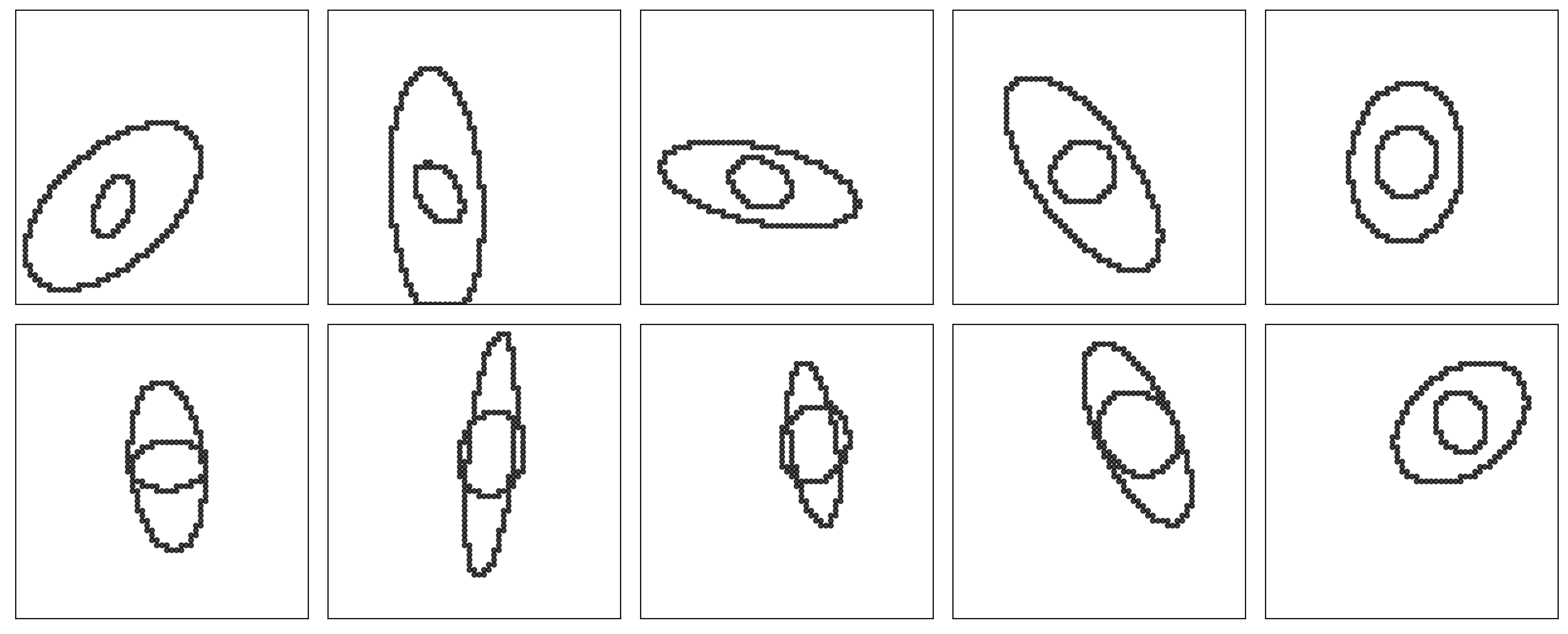}
		\caption{Ten measures from the nested-ellipse dataset.}
		\label{fig:ten_ellipses}
	\end{figure}

	We first consider the nested-ellipse benchmark used in \cite{cuturi2014fast,janati2020debiased,altschuler2021}, for which a reference optimal barycenter is available \citep{altschuler2021}. The dataset contains ten probability measures, each represented on a $60\times60$ grid and corresponding to one nested ellipse; see Figure~\ref{fig:ten_ellipses}. The marginal supports contain between $139$ and $192$ equally weighted atoms. The exact barycenter need not to be on the original grid. A direct MOT formulation can involve on the order of $192^{10}\approx10^{22}$ tuple variables and is therefore computationally intractable.
	
\begin{figure}[htbp]
		\centering
		\captionsetup[subfigure]{justification=centering}
		\begin{subfigure}[t]{0.24\textwidth}
			\centering
			\includegraphics[width=\textwidth]{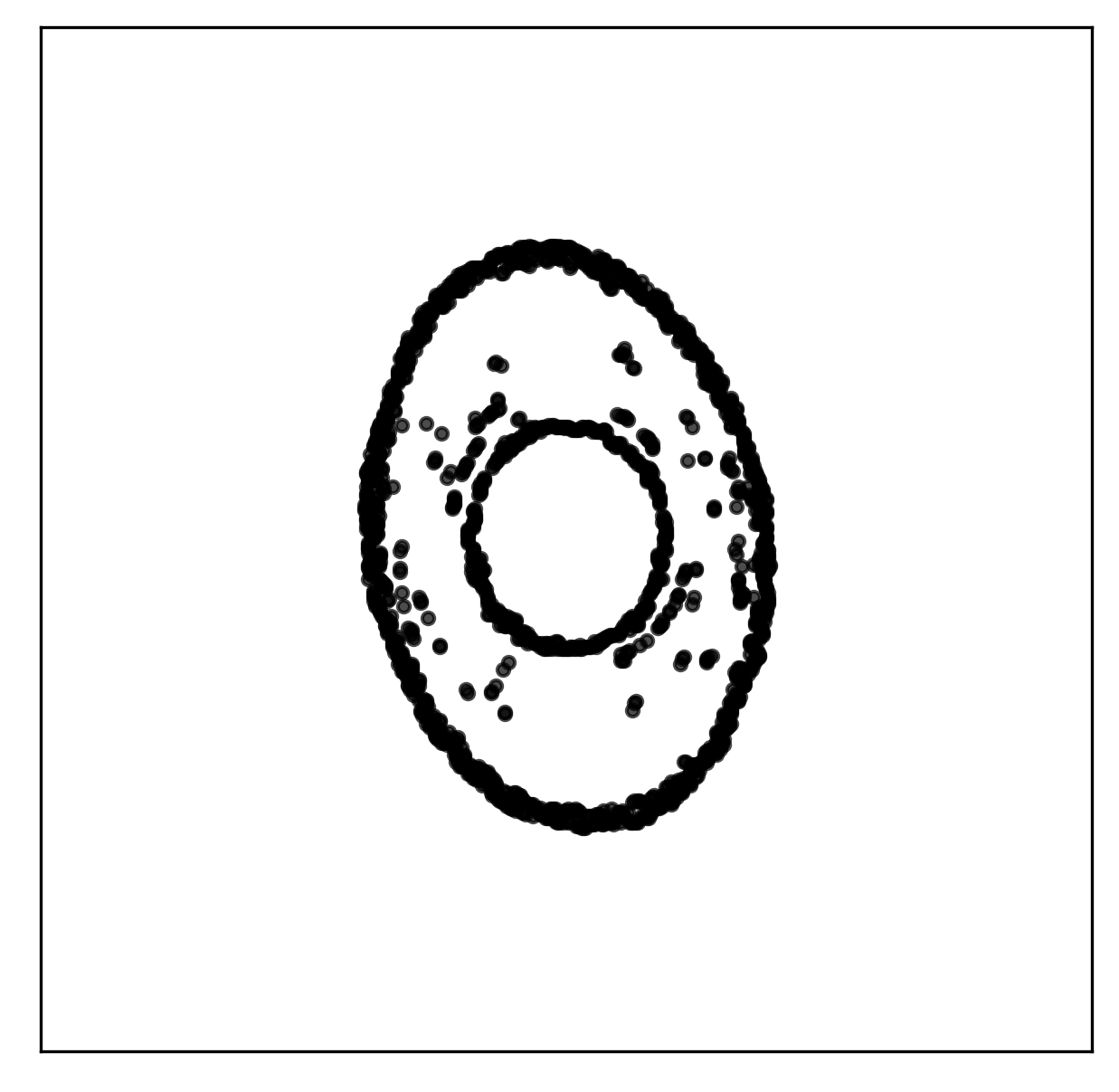}
			\caption{Exact ellipse \citep{altschuler2021}, objective=0.02666, time=8865.03s.}
		\end{subfigure}
		\hfill
		\begin{subfigure}[t]{0.24\textwidth}
			\centering
			\includegraphics[width=\textwidth]{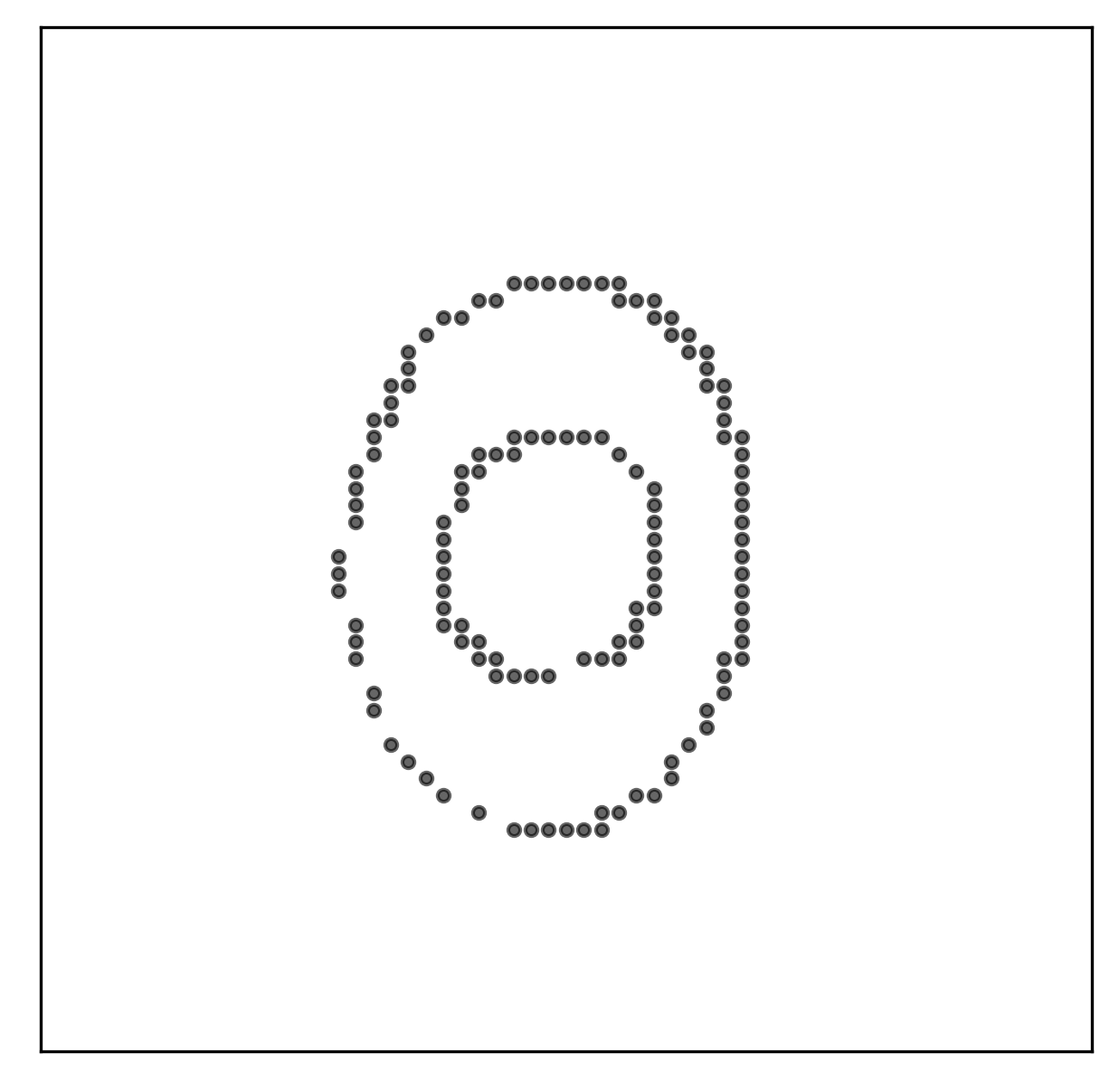}
			\caption{$10$ runs of Algorithm \ref{alg:subset_sampling} with $t=1$, objective=0.02709, time=17.46s.}
		\end{subfigure}
		\hfill
		\begin{subfigure}[t]{0.24\textwidth}
			\centering
			\includegraphics[width=\textwidth]{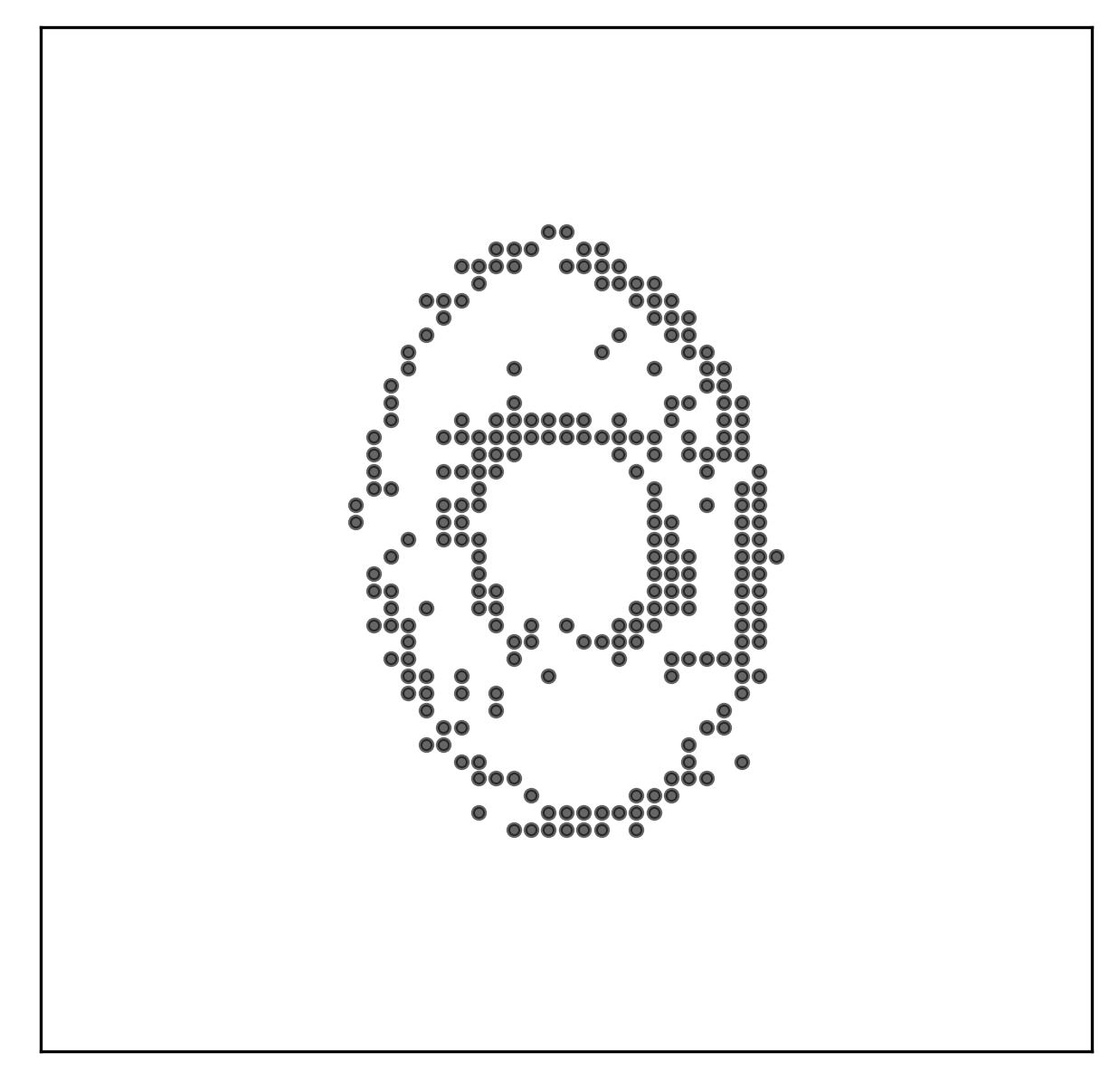}
			\caption{Algorithm \ref{alg:subset_enumeration} with $t=1$, objective=0.02673, time=21.08s.}
		\end{subfigure}
		\hfill
		\begin{subfigure}[t]{0.24\textwidth}
			\centering
			\includegraphics[width=\textwidth]{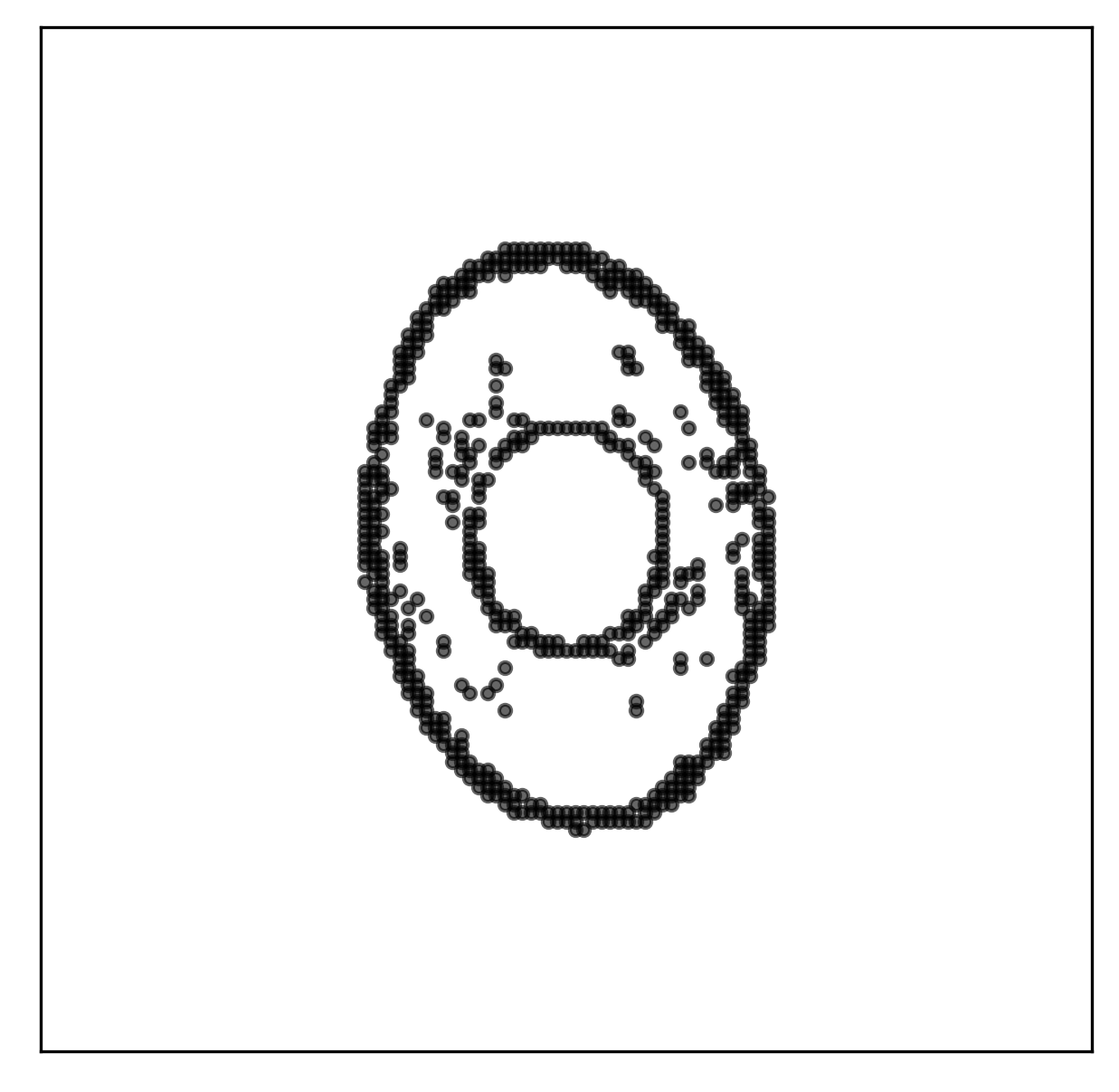}
			\caption{Algorithm \ref{alg:subset_enumeration} with $t=2$, objective=0.02667, time=212.44s.}
		\end{subfigure}
		\caption{Barycenters produced by different algorithms on the nested-ellipse dataset.}
		\label{fig:barycenters_nested_ellipse}
	\end{figure}

	We run $10$ runs of randomized sampling Algorithm \ref{alg:subset_sampling} and pick the best solution. We also use deterministic subset enumeration Algorithm \ref{alg:subset_enumeration}. Note that the deterministic candidate 
	support is the union of the candidate supports generated by all possible 
	sampling outcomes. It therefore contains the candidate support associated with 
	every individual realization and yields an objective value no larger than that 
	obtained from any single realization. Figure~\ref{fig:barycenters_nested_ellipse} illustrates this distinction: even the best of the ten possible $t=1$ samples is visibly worse than the deterministic $t=1$ solution. Accordingly, the remaining experiments focus on deterministic enumeration with $t=1$ and $t=2$.
	
	The empirical performance is substantially better than the worst-case guarantees in Theorem \ref{thm:W2_enumeration}. Moreover, small improvements in the objective can correspond to meaningful visual improvements: although $t=1$ is faster, $t=2$ produces a barycenter that is nearly indistinguishable from the reference solution.
	
	\subsection{MNIST dataset}
	To assess scalability in the number of marginals, we use the MNIST dataset \citep{lecun2002gradient}. For each of three digit classes, we randomly select $k=50$ empirical measures on a $28\times28$ grid ($d=2$). Their mass distributions differ, and the number of positive-mass pixels ranges from $165$ to $263$. The full MOT formulation can contain on the order of $263^{50}\approx10^{120}$ tuple variables.
	
\begin{figure}[htbp]
		\centering
		\captionsetup[subfigure]{justification=centering}
		\begin{subfigure}[b]{0.32\textwidth}
			\centering
			\includegraphics[width=\textwidth]{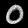}
			\caption{Algorithm \ref{alg:subset_enumeration_with_rep} with $t=1$, objective=0.0132, time=556s.}
		\end{subfigure}
		\hfill
		\begin{subfigure}[b]{0.32\textwidth}
			\centering
			\includegraphics[width=\textwidth]{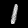}
			\caption{Algorithm \ref{alg:subset_enumeration_with_rep} with $t=1$, objective=0.0174, time=56s.}
		\end{subfigure}
		\hfill
		\begin{subfigure}[b]{0.32\textwidth}
			\centering
			\includegraphics[width=\textwidth]{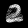}
			\caption{Algorithm \ref{alg:subset_enumeration_with_rep} with $t=1$, objective=0.0279, time=253s.}\end{subfigure}
		\vskip\baselineskip
		\begin{subfigure}[b]{0.32\textwidth}
			\centering
			\includegraphics[width=\textwidth]{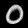}
			\caption{Algorithm \ref{alg:subset_enumeration_with_rep} with $t=2$, objective=0.0127, time=9213s.}
		\end{subfigure}
		\hfill
		\begin{subfigure}[b]{0.32\textwidth}
			\centering
			\includegraphics[width=\textwidth]{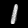}
			\caption{Algorithm \ref{alg:subset_enumeration_with_rep} with $t=2$, objective=0.0169, time=1148s.}
		\end{subfigure}
		\hfill
		\begin{subfigure}[b]{0.32\textwidth}
			\centering
			\includegraphics[width=\textwidth]{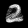}
			\caption{Algorithm \ref{alg:subset_enumeration_with_rep} with $t=2$, objective=0.0274, time=2716s.}
		\end{subfigure}
		\caption{Barycenters produced by different algorithms for MNIST.}
		\label{fig:barycenters_MNIST}
	\end{figure}

	We run Algorithm~\ref{alg:subset_enumeration_with_rep} with $t=1$ and $t=2$ and report the objective values and running times in Figure~\ref{fig:barycenters_MNIST}. Both settings are computationally viable, while $t=2$ consistently improves the objective value and preserves the characteristic shape of each digit more clearly.

	\subsection{Sign language dataset}
	The final experiment studies measures with large supports. Each $64\times64$ RGB image is represented as a five-dimensional empirical measure: two coordinates encode spatial location, and three encode color. For each of three hand gestures, we randomly select three training images as the input measures.


	\paragraph{Column generation implementation. }
	Since every image has a large support, explicitly constructing the full candidate support and solving the resulting LP in Algorithm~\ref{alg:subset_enumeration} is computationally expensive. 
	
	For $t=1$, the candidate support reduces to the union of the 
	marginal supports:
	\[
	S_2^1
	=
	\bigcup_{i\in[k]}\Xi_i.
	\]
	Thus, $S_2^1$ contains at most 
	$\sum_{i\in[k]}|\Xi_i|$ candidate atoms. For $t=2$, the candidate support 
	contains all pairwise centroids:
	\[
	S_{2}^{2}
	=
	\left\{
	\frac{\bm{x}_r+\bm{x}_s}{2}:
	\bm{x}_r\in\Xi_r,\ 
	\bm{x}_s\in\Xi_s,\ 
	1\le r<s\le k
	\right\}.
	\]
	Thus,
	\[
	|S_2^2|
	\leq
	\sum_{1\leq r<s\leq k}
	|\Xi_r|\,|\Xi_s|,
	\]
	so the size of $S_2^2$ grows quadratically with the marginal 
	support sizes. To control memory usage, we combine 
	Algorithm~\ref{alg:subset_enumeration} with a column-generation procedure 
	in the spirit of \cite{borgwardt2022column}. The restricted master problem 
	retains only a small set of active columns, and the pricing problem adds a 
	column only when its reduced cost is negative.
	Let 
	\[
	S_{\mathrm{cand}}\subseteq S_1^1\cup S_2^2
	\]
	be a mixed candidate support containing both individual input atoms and 
	pairwise centroids. A column is indexed by a candidate barycenter atom 
	$\bm{w}\in S_{\mathrm{cand}}$ and a tuple 
	$\vec{\bm{j}}=(j_1,\ldots,j_k)$ with 
	$j_i\in[|\Xi_i|]$. Its cost is
	\[
	c(\bm{w},\vec{\bm{j}})
	=
	\sum_{i\in[k]}
	\lambda_i
	\left\|
	\bm{w}-\hat{\bm{x}}_{i j_i}
	\right\|_2^2 .
	\]
	
	For a current column set $\mathcal C$, the restricted master problem is
	\[
	\begin{aligned}
		\min_{\bm\theta\ge 0}\quad
		&
		\sum_{(\bm w,\vec{\bm{j}})\in\mathcal C}
		c(\bm{w},\vec{\bm{j}})\theta_{\bm{w},\vec{\bm{j}}}
		\\
		\mathrm{s.t.}\quad
		&
		\sum_{\substack{(\bm{w},\vec{\bm{j}})\in\mathcal C: j_i=j}}
		\theta_{\bm{w},\vec{\bm{j}}}
		=
		\hat p_{ij},
		&&
		\forall i\in[k],\ j\in[|\Xi_i|].
	\end{aligned}
	\tag{CG-RMP}
	\]
	The variable $\theta_{\bm w,\vec{\bm j}}$ is the mass assigned to barycenter atom $\bm w$ through the tuple $(\hat{\bm x}_{1j_1},\ldots,\hat{\bm x}_{kj_k})$. The induced mass of atom $\bm w$ and the corresponding barycenter are
	\[
	m_{\bm{w}}
	=
	\sum_{\vec{\bm{j}}:(\bm{w},\vec{\bm{j}})\in\mathcal C}
	\theta_{\bm{w},\vec{\bm{j}}}, \quad
	\mathbb{P}_{\mathcal C}
	=
	\sum_{\bm{w}\in S_{\mathrm{cand}}:m_{\bm{w}}>0} m_{\bm{w}}\delta_{\bm{w}} .
	\]
	Let $\gamma_{ij}$ be the dual variable for the marginal constraint associated with $\hat{\bm x}_{ij}$. The reduced cost of $(\bm w,\vec{\bm j})$ is
	\[
	\bar c(\bm{w},\vec{\bm{j}})
	=
	\sum_{i\in[k]}
	\lambda_i
	\left\|
	\bm{w}-\hat{\bm{x}}_{i j_i}
	\right\|_2^2
	-
	\sum_{i\in[k]}
	\gamma_{i j_i}.
	\]
	For a fixed $\bm w$, the minimizing tuple decomposes across marginals:
	\[
	j_i(\bm{w})
	\in
	\arg\min_{j\in[|\Xi_i|]}
	\left\{
	\lambda_i
	\left\|
	\bm{w}-\hat{\bm{x}}_{ij}
	\right\|_2^2
	-
	\gamma_{ij}
	\right\},
	\forall i\in[k].
	\]
	Hence, the pricing value for $\bm w$ is
	\[
	r(\bm{w})
	=
	\sum_{i\in[k]}
	\min_{j\in[|\Xi_i|]}
	\left\{
	\lambda_i
	\left\|
	\bm{w}-\hat{\bm{x}}_{ij}
	\right\|_2^2
	-
	\gamma_{ij}
	\right\}.
	\]
	If $\min_{\bm w\in S_{\mathrm{cand}}}r(\bm w)<0$, we add the corresponding column $(\bm w,j_1(\bm w),\ldots,j_k(\bm w))$ and reoptimize the restricted master problem. If no negative-reduced-cost column exists, the current solution is optimal for the LP restricted to $S_{\mathrm{cand}}$.

	For the sign-language instances, even scanning all of $S_2^2$ is expensive. We therefore construct a hybrid support intended to optimize the accuracy of $t=2$ at a computational cost closer to $t=1$. First, we solve the $t=1$ restricted problem and collect its active atoms $A_1=\{\bm w:m_{\bm w}>0\}$. For each $\bm a\in A_1$ and marginal $i\in[k]$, let $\mathcal N_i^K(\bm a)\subseteq\Xi_i$ contain the $K$ nearest atoms of $\mathbb P_i$ to $\bm a$; we set $K=5$. Define
	\[
	S_{\mathrm{hyb}}
	=
	A_1
	\cup
	\bigcup_{\bm{a}\in A_1}
	\left\{
	\frac{\bm{x}_r+\bm{x}_s}{2}:
	\bm{x}_r\in\mathcal N_r^K(\bm{a}),\
	\bm{x}_s\in\mathcal N_s^K(\bm{a}),\
	\forall 1\le r<s\le k
	\right\}.
	\]
	The pricing step is performed over
	$S_{\mathrm{cand}}=S_{\mathrm{hyb}}$ rather than the full candidate set
	$S_2^2$. The hybrid set contains the active atoms of the $t=1$ solution
	and pairwise centroids formed from their neighboring input atoms. This localized construction substantially reduces the number of candidate
	atoms while retaining potentially useful $t=2$ support points.
	
	We warm-start column generation with the active $t=1$ columns, which provide a feasible restricted master problem. The procedure terminates when either no negative-reduced-cost column is found or the time limit is reached. Under a time limit, the final restricted-master solution remains a feasible barycenter, and we report its objective value. The worst-case guarantee of Algorithm~\ref{alg:subset_enumeration} requires optimizing over the full enumerated support; it does not apply to this hybrid acceleration. We refer to the implementation as ``Hybrid Algorithm~\ref{alg:subset_enumeration} with $t=2$.''

\begin{figure}[htbp]
		\centering
		\captionsetup[subfigure]{justification=centering}
		\begin{subfigure}[t]{0.24\textwidth}
			\centering
			\includegraphics[width=\textwidth]{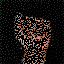}
			\caption{Algorithm \ref{alg:subset_enumeration} with $t=1$, objective=0.00368, time=1h.}
			\label{fig:barycenter_sign_A_(a)}
		\end{subfigure}
		\hfill
		\begin{subfigure}[t]{0.24\textwidth}
			\centering
			\includegraphics[width=\textwidth]{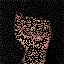}
			\caption{Algorithm \ref{alg:subset_enumeration} with $t=1$, objective=0.00358, time=2h.}
			\label{fig:barycenter_sign_A_(b)}
		\end{subfigure}
		\hfill
		\begin{subfigure}[t]{0.24\textwidth}
			\centering
			\includegraphics[width=\textwidth]{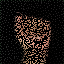}
			\caption{Hybrid Algorithm \ref{alg:subset_enumeration} with $t=2$, objective=0.00355, time=1h.}
			\label{fig:barycenter_sign_A_(c)}
		\end{subfigure}
		\hfill
		\begin{subfigure}[t]{0.24\textwidth}
			\centering
			\includegraphics[width=\textwidth]{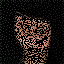}
			\caption{Hybrid Algorithm \ref{alg:subset_enumeration} with $t=2$, objective=0.00348, time=2h.}
			\label{fig:barycenter_sign_A_(d)}
		\end{subfigure}
		\caption{Barycenters produced by different algorithms for gesture 1.}
		\label{fig:barycenter_sign_A}
	\end{figure}

\begin{figure}[htbp]
		\centering
		\captionsetup[subfigure]{justification=centering}
		\begin{subfigure}[t]{0.24\textwidth}
			\centering
			\includegraphics[width=\textwidth]{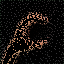}
			\caption{Algorithm \ref{alg:subset_enumeration} with $t=1$, objective=0.0115, time=1h.}
		\end{subfigure}
		\hfill
		\begin{subfigure}[t]{0.24\textwidth}
			\centering
			\includegraphics[width=\textwidth]{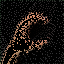}
			\caption{Algorithm \ref{alg:subset_enumeration} with $t=1$, objective=0.0115, time=2h.}
		\end{subfigure}
		\hfill
		\begin{subfigure}[t]{0.24\textwidth}
			\centering
			\includegraphics[width=\textwidth]{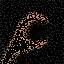}
			\caption{Hybrid Algorithm \ref{alg:subset_enumeration} with $t=2$, objective=0.0111, time=1h.}
		\end{subfigure}
		\hfill
		\begin{subfigure}[t]{0.24\textwidth}
			\centering
			\includegraphics[width=\textwidth]{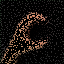}
			\caption{Hybrid Algorithm \ref{alg:subset_enumeration} with $t=2$, objective=0.0110, time=2h.}
		\end{subfigure}
		\caption{Barycenters produced by different algorithms for gesture 2.}
		\label{fig:barycenter_sign_C}
	\end{figure}

\begin{figure}[htbp]
		\centering
		\captionsetup[subfigure]{justification=centering}
		\begin{subfigure}[t]{0.24\textwidth}
			\centering
			\includegraphics[width=\textwidth]{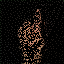}
			\caption{Algorithm \ref{alg:subset_enumeration} with $t=1$, objective=0.00957, time=1h.}
		\end{subfigure}
		\hfill
		\begin{subfigure}[t]{0.24\textwidth}
			\centering
			\includegraphics[width=\textwidth]{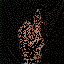}
			\caption{Algorithm \ref{alg:subset_enumeration} with $t=1$, objective=0.00949, time=2h.}
		\end{subfigure}
		\hfill
		\begin{subfigure}[t]{0.24\textwidth}
			\centering
			\includegraphics[width=\textwidth]{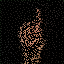}
			\caption{Hybrid Algorithm \ref{alg:subset_enumeration} with $t=2$, objective=0.00917, time=1h.}
		\end{subfigure}
		\hfill
		\begin{subfigure}[t]{0.24\textwidth}
			\centering
			\includegraphics[width=\textwidth]{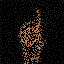}
			\caption{Hybrid Algorithm \ref{alg:subset_enumeration} with $t=2$, objective=0.00898, time=2h.}
		\end{subfigure}
		\caption{Barycenters produced by different algorithms for gesture 3.}
		\label{fig:barycenter_sign_D}
	\end{figure}

	\paragraph{Results.}
	Figures~\ref{fig:barycenter_sign_A}--\ref{fig:barycenter_sign_D}
	show representative results for the three gestures. Relative to the best $t=1$ solutions reported in the figures, the best hybrid $t=2$ solutions
	reduce the objective values by approximately $2.8\%$, $4.3\%$, and $5.4\%$
	for gestures~1, 2, and 3, respectively. The hybrid $t=2$ barycenters also
	provide visually clearer representations of the corresponding hand shapes.
	For gesture~1, Figures~\ref{fig:barycenter_sign_A_(a)} and
	\ref{fig:barycenter_sign_A_(b)} show the $t=1$ solutions obtained after one
	and two hours, respectively. Figures~\ref{fig:barycenter_sign_A_(c)} and
	\ref{fig:barycenter_sign_A_(d)} show the hybrid $t=2$ solutions obtained
	after one and two additional hours, respectively, using the one-hour $t=1$
	solution as a warm start. The hybrid $t=2$ barycenters form fuller and more
	clearly delineated hand shapes than their $t=1$ counterparts.
	Figures~\ref{fig:barycenter_sign_C} and \ref{fig:barycenter_sign_D} exhibit
	similar improvements for gestures~2 and~3.
	
	We also evaluate the barycenters in a nearest-barycenter classification task. Each gesture class is represented by its training-sample barycenter, and each of the $25$ test images per gesture is assigned to the class with the smallest Wasserstein distance. The two-hour $t=1$ barycenters yield $96\%$ accuracy, whereas the one-hour hybrid $t=2$ barycenters classify all $75$ images correctly. Thus, the improvement in the barycenter objective also translates into better downstream classification performance; see \Cref{table:gesture_class}.
	
	\begin{table}[htbp]
		\centering\small
		\caption{Gesture-classification results. Each entry is the number of test images in the indicated category. For example, Algorithm~\ref{alg:subset_enumeration} with $t=1$ correctly classifies $22$ of the $25$ gesture-2 images and misclassifies the remaining $3$.}
		\label{table:gesture_class}
		\begin{tabular}{l|r|r|r}
			\hline
			Methods       & gesture 1 & gesture 2 & gesture 3 \\ \hline
			Algorithm \ref{alg:subset_enumeration} with $t=1$ correct & 25        & 22        & 25        \\ \hline
			Algorithm \ref{alg:subset_enumeration} with $t=1$ wrong   & 0         & 3         & 0         \\ \hline
			Hybrid Algorithm \ref{alg:subset_enumeration} with $t=2$ correct & 25        & 25        & 25        \\ \hline
			Hybrid Algorithm \ref{alg:subset_enumeration} with $t=2$ wrong   & 0         & 0         & 0         \\ \hline
		\end{tabular}
	\end{table}

	\section{Extensions}\label{section:extension}
	
	We consider two extensions of the candidate-support framework. First, we show that the deterministic support $S_1^t$ yields a value guarantee for sparse type-$2$ barycenters, although the resulting restricted problem remains a mixed-integer program. Second, we exploit coordinatewise weighted medians to construct exact and approximate supports for type-$1$ barycenters. Extending comparable guarantees to general type-$p$ barycenters with $p\notin\{1,2\}$ remains open.
	
	\subsection{Approximation guarantee for sparse type-$2$ Wasserstein barycenters}
	\label{subsec:sparse-w2-support-reduction}

	Consider the support-size-constrained Wasserstein barycenter problem:
	\begin{align*}
		v_{m}^{\star}:=\inf_{\substack{\mathbb{P}\in\mathcal P(\mathbb R^d)\\
				|\text{supp}(\mathbb{P})|\le m}}
		\sum_{i\in[k]}\lambda_i W_2^2(\mathbb{P},\mathbb{P}_i).
	\end{align*}
	For a finite candidate set $G\subseteq\mathbb R^d$, define the value of the corresponding restricted sparse Wasserstein barycenter problem as
	\begin{align*}
		v_{m}(G):=\inf_{\substack{\mathbb{P}\in\mathcal P(\mathbb R^d)\\
				\text{supp}(\mathbb{P})\subseteq G,\,
				|\text{supp}(\mathbb{P})|\le m}}
		\sum_{i\in[k]}\lambda_i W_2^2(\mathbb{P},\mathbb{P}_i).
	\end{align*}
	\paragraph{Mixed-integer formulation of the restricted problem.}
	\label{subsec:restricted-sparse-formulation}
	Let $S=\{\bm s_1,\ldots,\bm s_M\}\subseteq\mathbb R^d$ be a finite 
	candidate support, and write
	$\mathbb P_i=\sum_{j\in[|\Xi_i|]}\hat p_{ij}
	\delta_{\hat{\bm x}_{ij}}$ for each $i\in[k]$.
	For each candidate atom $\bm s_\ell$, let $z_\ell$ denote its barycenter 
	mass and let $y_\ell\in\{0,1\}$ be 
	a binary activation variable for that atom.
	For each $i\in[k]$, $j\in[|\Xi_i|]$, and $\ell\in[M]$, let $\pi_{ij\ell}$ denote the amount of mass transported from the 
	barycenter atom $\bm s_\ell$ to the input atom $\hat{\bm x}_{ij}$ in the 
	transport plan between the barycenter and $\mathbb P_i$.
	The restricted $m$-sparse problem over $S$ can then be formulated as the following 
	mixed-integer linear program:
	\[
	\begin{aligned}
		v_{m}(S)
		=
		\min_{\bm\pi,\bm{z},\bm{y}}\quad
		&
		\sum_{i\in[k]}\lambda_i
		\sum_{j\in[|\Xi_i|]}
		\sum_{\ell\in[M]}
		\left\|\bm{s}_\ell-\hat{\bm{x}}_{ij}\right\|_2^2
		\pi_{ij\ell}
		\\
		\mathrm{s.t.}\quad
		&
		\sum_{\ell\in[M]}\pi_{ij\ell}
		=
		\hat p_{ij},
		&&
		\forall i\in[k],\ j\in[|\Xi_i|],
		\\
		&
		\sum_{j\in[|\Xi_i|]}\pi_{ij\ell}
		=
		z_\ell,
		&&
		\forall i\in[k],\ \ell\in[M],
		\\
		&
		0\le z_\ell\le y_\ell,
		&&
		\forall \ell\in[M],
		\\
		&
		\sum_{\ell\in[M]}y_\ell\le m,
		\\
		&
		\pi_{ij\ell}\ge0,\quad y_\ell\in\{0,1\},
		&&
		\forall i\in[k],\ j\in[|\Xi_i|],\ \ell\in[M].
	\end{aligned}
	\tag{RWBCenter}
	\label{eq:restricted-sparse-barycenter}
	\]
	The first constraint family ensures that the transport 
	plan associated with each marginal $\mathbb P_i$ has the prescribed 
	marginal masses $\hat p_{ij}$. The second constraint family requires all 
	$k$ transport plans to share the same barycenter mass vector 
	$\bm z=(z_1,\ldots,z_M)$.
	The linking constraints $z_\ell\leq y_\ell$ deactivate unselected atoms, 
	and $\sum_{\ell\in[M]}y_\ell\leq m$ limits the support size. The final constraints impose nonnegativity and 
	integrality.
	
	Given an optimal solution 
	$(\bm\pi^\star,\bm z^\star,\bm y^\star)$, the corresponding restricted 
	sparse barycenter is
	\[
	\mathbb P_{S,m}
	=
	\sum_{\ell\in[M]}z_\ell^\star\delta_{\bm s_\ell}.
	\]
	Candidate atoms with $z_\ell^\star=0$ can be omitted from 
	its support.
	We next analyze the value obtained when this MILP is solved to optimality 
	over the candidate support $S_1^t$.
	
	\paragraph{Approximation guarantee.}
	We next show that the deterministic candidate support
	$S_1^t$ defined in
	\eqref{eq:candidate_support_with_rep_deterministic} can also be used to
	approximate sparse barycenters. The following lemma provides the key
	ingredient: it shows that $S_1^t$ yields a
	$(1+1/t)$-approximation for the weighted one-center problem induced by
	any probability distribution over the indexed input atoms.
	\begin{lemma}
		\label{lem:w2-centroid-set}
		For every matrix $\bm{q}\in \mathbb{R}_+^{k \times n}$ satisfying $\sum_{i\in[k]}\sum_{j\in \Xi_i}q_{ij}=1$, we have
		\begin{align*}
			\min_{\bm{w} \in S^{t}_1}
			\sum_{i \in [k]}\sum_{j \in \Xi_i}
			q_{ij}\|\bm{w}-\hat{\bm{x}}_{ij}\|_2^2
			\leq
			\left(1+\frac1t\right)
			\min_{\bm{w}\in\mathbb R^d}
			\sum_{i \in [k]}\sum_{j \in \Xi_i}
			q_{ij}\|\bm{w}-\hat{\bm{x}}_{ij}\|_2^2.
		\end{align*}
	\end{lemma}
	\begin{proof}
		Let $\bm{Y}$ be a random variable supported on
		$\{\hat{\bm{x}}_{ij}:i\in [k],j \in \Xi_i\}$ with $\mathbb P(\bm{Y}=\hat{\bm{x}}_{ij})=q_{ij}$. Let $\bm{\mu}:=\mathbb E[\bm{Y}]$. Then $\bm{\mu} \in \arg\min_{\bm{w}\in\mathbb R^d} \mathbb E\|\bm{w}-\bm{Y}\|_2^2$,
		and $\min_{\bm{w}\in\mathbb R^d} \mathbb E\|\bm{w}-\bm{Y}\|_2^2=\mathbb E\|\bm{\mu}-\bm{Y}\|_2^2$. Now draw $\bm{Y}_1,\ldots,\bm{Y}_t$ independently from the distribution of $\bm{Y}$, and define
		$
		\bar{\bm{Y}}_t:=\frac1t\sum_{r=1}^t \bm{Y}_r.$
		By taking expectation over $\bm{Y},\bm{Y}_1,\ldots,\bm{Y}_t$, and expand the squares similar to \eqref{eq:l2_expansion}, the approximation follows from the proof of Theorem \ref{thm:W2_sampling_with_rep}.
	\end{proof}
	We next transfer this local centroid guarantee to the sparse barycenter objective.
	\begin{thm}
		\label{thm:sparse-w2-gt}
		For every integer $t\ge1$,
		\[
		v_{m}(S_1^t)
		\le
		\left(1+\frac1t\right)v_{m}^{\star}.
		\]
	\end{thm}
	
	\begin{proof}
		Choose an optimal $m$-sparse barycenter $\mathbb{P}=\sum_{\ell=1}^m a_\ell\delta_{\bm{w}_\ell}$ such that $
		\sum_{i\in[k]}\lambda_i W_2^2(\mathbb{P},\mathbb{P}_i)=v_{m}^{\star}$ and $\bm{a}\in \mathbb{R}_{++}^m,\sum_{\ell=1}^m a_\ell=1$.
		For each $i\in[k]$, let $
		\bm{\pi}_i=\{\pi_{i\ell j}\}_{\ell\in[m],\,j\in[|\Xi_i|]}$
		be an optimal transport plan between $\mathbb{P}$ and $\mathbb{P}_i$. Hence,
		\begin{align*}
			\sum_{j\in[|\Xi_i|]}\pi_{i\ell j}=a_\ell,
			\qquad \sum_{\ell=1}^m\pi_{i\ell j}=\hat p_{ij}.
		\end{align*}
		
		For each active barycenter atom $\bm{w}_\ell$, define the probability mass
		$
		q^\ell_{ij}:=\frac{\lambda_i\pi_{i\ell j}}{a_\ell}$ for each $i\in [k]$ and $j\in [|\Xi_i|]$. 
		This is indeed a probability vector, since
		\begin{align*}
			\sum_{i \in [k]}\sum_{j \in [|\Xi_i|]}q^\ell_{ij}=
			\frac{1}{a_\ell}
			\sum_{i\in[k]}\lambda_i
			\sum_{j\in[|\Xi_i|]}\pi_{i\ell j}=\frac{1}{a_\ell}
			\sum_{i\in[k]}\lambda_i a_\ell=1.
		\end{align*}
		
		By Lemma~\ref{lem:w2-centroid-set}, for each $\ell\in[m]$, there exists $\bm{s}_\ell\in S_1^t$ such that
		\begin{align*}
			\sum_{i \in [k]}\sum_{j \in [|\Xi_i|]}
			q^\ell_{ij}\|\bm{s}_\ell-\hat{\bm{x}}_{ij}\|_2^2
			\leq
			\left(1+\frac1t\right)
			\min_{\bm{w}\in\mathbb R^d}
			\sum_{i \in [k]}\sum_{j \in [|\Xi_i|]}
			q^\ell_{ij}\|\bm{w}-\hat{\bm{x}}_{ij}\|_2^2.
		\end{align*}
		In particular, since $\bm{w}_\ell\in\mathbb R^d$ is feasible to the right-hand minimization problem, we have
		\begin{align} \label{eq:sparse_local_approx}
			\sum_{i \in [k]}\sum_{j \in [|\Xi_i|]}
			q^\ell_{ij}\|\bm{s}_\ell-\hat{\bm{x}}_{ij}\|_2^2
			\le
			\left(1+\frac1t\right)
			\sum_{i \in [k]}\sum_{j \in [|\Xi_i|]}
			q^\ell_{ij}\|\bm{w}_\ell-\hat{\bm{x}}_{ij}\|_2^2.
		\end{align}
		
		Now define $
		\tilde{\mathbb{P}}:=\sum_{\ell=1}^ma_\ell\delta_{\bm{s}_\ell}.$
		If several $\bm{s}_\ell$'s coincide, we merge their masses.  Therefore,
		$\tilde{\mathbb{P}}$ is supported on $S_1^t$ and has support size at most
		$m$.  We construct a feasible transportation plan from $\tilde{\mathbb{P}}$
		to each $\mathbb{P}_i$ by transporting the mass $a_\ell$ at $\bm{s}_\ell$ to the
		points $\hat{\bm{x}}_{ij}$ using the same coefficients $\pi_{i \ell j}$.
		Thus,
		\begin{align*}
			v_{m}(S_1^t)
			&\leq
			\sum_{i\in[k]}\lambda_i W_2^2(\tilde{\mathbb{P}},\mathbb{P}_i)
			\leq
			\sum_{i\in[k]}\lambda_i
			\sum_{\ell=1}^m
			\sum_{j\in[|\Xi_i|]}
			\pi_{i\ell j}
			\|\bm{s}_\ell-\hat{\bm{x}}_{ij}\|_2^2
			\\
			&=
			\sum_{\ell=1}^m
			a_\ell
			\sum_{i \in [k]}\sum_{j \in \Xi_i}
			q^\ell_{ij}
			\|\bm{s}_\ell-\hat{\bm{x}}_{ij}\|_2^2
			\leq
			\left(1+\frac1t\right)
			\sum_{\ell=1}^m
			a_\ell
			\sum_{i \in [k]}\sum_{j \in \Xi_i}
			q^\ell_{ij}
			\|\bm{w}_\ell-\hat{\bm{x}}_{ij}\|_2^2
			\\
			&=
			\left(1+\frac1t\right)
			\sum_{i\in[k]}\lambda_i
			\sum_{\ell=1}^m
			\sum_{j\in[|\Xi_i|]}
			\pi_{i\ell j}
			\|\bm{w}_\ell-\hat{\bm{x}}_{ij}\|_2^2
			\\
			&=
			\left(1+\frac1t\right)
			\sum_{i\in[k]}\lambda_i W_2^2(\mathbb{P},\mathbb{P}_i)
			=\left(1+\frac{1}{t}\right)
			v_{m}^{\star},
		\end{align*}
		where the first inequality is since $\supp(\tilde{\mathbb{P}}) \subseteq S_1^t$, the second inequality is since $\bm{\pi}_i$ is a feasible transportation plan from $\tilde{\mathbb{P}}$ to $\mathbb{P}_i$, the first and second equalities are by the definition of $q_{ij}^{\ell}$, the third inequality follows \eqref{eq:sparse_local_approx}, the third equality is by the definition of $(\pi_{i \ell j},\bm{w}_\ell)$, and the last inequality is by the optimality of $\mathbb{P}$. This completes the proof.
	\end{proof}
	Theorem~\ref{thm:sparse-w2-gt} shows that restricting the barycenter support to $S_1^t$ increases the optimal $m$-sparse barycenter value by at most a factor of $1+1/t$. This result should be interpreted as a candidate-support reduction guarantee rather than a PTAS for the sparse barycenter problem. Indeed, unlike the nonsparse problem studied earlier,
	restricting the candidate support does not reduce the sparse problem to an
	LP: computing $v_m(S_1^t)$ still requires solving the MILP
	\eqref{eq:restricted-sparse-barycenter} to optimality. Nevertheless, the
	reduced candidate support can substantially decrease the size of the
	resulting formulation and therefore improve its practical tractability.

	\subsection{Approximation guarantee for the type-$1$ Wasserstein barycenter problem}
	
	The discrete type-$1$ Wasserstein barycenter problem is
	\begin{align*}
		v_I^*
		=
		\inf_{\mathbb{P}\in\mathcal P_1(\mathbb R^d)}
		\sum_{i\in[k]}
		\lambda_i W_1(\mathbb P,\mathbb P_i),
	\end{align*}
	where $\mathcal P_1(\mathbb R^d)$ is the set of Borel probability
	measures with finite first moment and
	\begin{align*}
		W_1(\mathbb P,\mathbb Q)
		:=
		\inf_{\Pi\in\mathcal M(\mathbb P,\mathbb Q)}
		\int_{\mathbb R^d\times\mathbb R^d}
		\|\bm x-\bm y\|_1\,d\Pi(\bm x,\bm y).
	\end{align*}
	
	The MOT representation underlying
	\Cref{lem:optimal_support} extends directly to the $\ell_1$
	transportation cost. For each tuple
	$(\bm x_1,\ldots,\bm x_k)$, its tuplewise transportation cost is
	\begin{align*}
		\min_{\bm s\in\mathbb R^d}
		\sum_{i\in[k]}
		\lambda_i\|\bm s-\bm x_i\|_1.
	\end{align*}
	Unlike the squared Euclidean case, the minimizer need not
	be unique. Nevertheless, the coordinatewise separability of the
	$\ell_1$ norm allows us to construct a finite exact candidate support. Indeed, we have
	\begin{align*}
		\min_{\bm s\in\mathbb R^d}
		\sum_{i\in[k]}\lambda_i\|\bm s-\bm x_i\|_1
		=
		\sum_{j\in[d]}
		\min_{s_j\in\mathbb R}
		\sum_{i\in[k]}\lambda_i|s_j-x_{ij}|.
	\end{align*}
	Thus, a tuplewise minimizer can be obtained by solving $d$ one-dimensional weighted-median problems
	\citep{vazler2012weighted,sabo2008best}.
	\begin{lemma}[Weighted median in one dimension; Lemma~3.1 in \citealt{sabo2008best}]
		\label{lemma:weighted_median}
		Let $a_1\leq a_2\leq\cdots\leq a_k$ be real numbers with corresponding weights $\lambda_1,\ldots,\lambda_k$, where $\sum_{i\in[k]}\lambda_i=1$. Define
		$
		J
		=
		\left\{
		\nu\in[k]:
		\sum_{i\in[\nu]}\lambda_i\leq \frac{1}{2}
		\right\}.
		$
		If $J\neq\emptyset$, let $\nu_0=\max_{\nu\in J}\nu$. Then:
		\begin{itemize}
			\item[(i)] If $J=\emptyset$, then
			$
			\min_{s\in\mathbb{R}}
			\sum_{i\in[k]}\lambda_i|a_i-s|$
			is attained at $s^*=a_1$.
			
			\item[(ii)] If $J\neq\emptyset$ and
			$
			\sum_{i=1}^{\nu_0}\lambda_i<\frac{1}{2},
			$
			then the minimum is attained at $s^*=a_{\nu_0+1}$.
			
			\item[(iii)] If $J\neq\emptyset$ and
			$
			\sum_{i=1}^{\nu_0}\lambda_i=\frac{1}{2},
			$
			then every point $s^*\in[a_{\nu_0},a_{\nu_0+1}]$ is optimal.
		\end{itemize}
	\end{lemma}
	For each coordinate $j\in[d]$, define the set of observed
	coordinate values
	\begin{align*}
		U_j
		:=
		\left\{
		x_j:
		\bm x\in\bigcup_{i\in[k]}\Xi_i
		\right\}.
	\end{align*}
	We define the coordinate-product candidate support
	\begin{align*}
		\hat S_1
		=
		\left\{
		\bm y\in\mathbb R^d:
		y_j\in U_j
		\text{ for every }j\in[d]
		\right\}.
	\end{align*}
	By \Cref{lemma:weighted_median}, for every input tuple
	$(\bm x_1,\ldots,\bm x_k)$, each coordinate of a tuplewise minimizer can
	be chosen from the corresponding coordinate values appearing in that
	tuple. Such values belong to $U_j$, and therefore $\hat S_1$ contains a
	tuplewise minimizer for every input tuple. Equivalently, for every
	$(\bm x_1,\ldots,\bm x_k)$, we have
	\begin{align*}
		\min_{\bm s\in\hat S_1}
		\sum_{i\in[k]}\lambda_i\|\bm s-\bm x_i\|_1
		=
		\min_{\bm s\in\mathbb R^d}
		\sum_{i\in[k]}\lambda_i\|\bm s-\bm x_i\|_1.
	\end{align*}
	Thus, restricting the type-$1$ MOT formulation
	to $\hat S_1$ leaves every tuplewise cost unchanged and therefore
	preserves the optimal barycenter value. Moreover, since
	\begin{align*}
		|U_j|
		\leq
		\sum_{i\in[k]}|\Xi_i|
		\leq nk,
		\qquad j\in[d],
	\end{align*}
	we have
	\begin{align*}
		|\hat S_1|
		=
		\prod_{j\in[d]}|U_j|
		\leq
		\left(
		\sum_{i\in[k]}|\Xi_i|
		\right)^d
		\leq
		(nk)^d.
	\end{align*}
	The coordinate-product structure of $\hat S_1$ motivates the reduced candidate supports introduced next.
	
	\paragraph{Subspace-sampling algorithm.} 
	Unlike the type-$2$ constructions, which sample marginals, the type-$1$ method samples coordinates. Let $\hat{U}:=\bigcup_{i\in[k]}\Xi_i$ and, for each $\ell\in[d]$, let $U_\ell:=\{x_\ell:\bm x\in \hat{U}\}$. For $\mathcal T\subseteq[d]$, write $\bm x_{\mathcal T}$ for the corresponding subvector. Fix $t\in\{0,\ldots,d\}$ and uniformly sample $\mathcal T\subseteq[d]$ with $|\mathcal T|=d-t$. Define
	
	\begin{align*}
		S_3^{\mathcal{T}}
		:=
		\left\{
		\bm{y}\in\mathbb{R}^d:
		\text{there exists }\bm{x}\in \hat{U}\text{ such that }
		\bm{y}_{\mathcal{T}}=\bm{x}_{\mathcal{T}},
		\text{ and }y_{\ell}\in U_{\ell}
		\text{ for every }\ell\notin\mathcal{T}
		\right\}.
	\end{align*}
	The coordinates in $\mathcal T$ are inherited jointly from one observed support point, whereas the remaining $t$ coordinates are selected independently from their observed values. Since $|\hat{U}|\leq nk$ and $|U_\ell|\leq nk$, the support size satisfies
	\begin{align*}
		|S_3^{\mathcal{T}}|
		\leq
		(nk)(nk)^t
		=
		(nk)^{t+1}.
	\end{align*}
	
	\begin{algorithm}[htbp]
		\caption{Subspace-sampling algorithm for the type-$1$ Wasserstein barycenter problem}
		\label{alg:subspace_sampling}
		\begin{algorithmic}[1]
			\State \textbf{Input}: Probability measures $\{\mathbb{P}_i\}_{i\in[k]}$ with supports $\{\Xi_i\}_{i\in[k]}$ and an integer $t\in\{0,\ldots,d\}$.
			\State Uniformly sample $\mathcal T\subseteq[d]$ with $|\mathcal T|=d-t$.
			\State Construct the candidate support $S_3^{\mathcal{T}}$.
			\State Solve the type-$1$ restricted MOT problem over $S_3^{\mathcal{T}}$ and recover an optimal primal solution $\bm{\Pi}^*(S_3^{\mathcal{T}})$.
			\State \textbf{Output}: Construct and return the approximate barycenter $\mathbb{P}(S_3^{\mathcal{T}})$ using \Cref{lem_distr}.
		\end{algorithmic}
	\end{algorithm}
	
	\begin{thm}\label{thm:W1_sampling}
		Let $S_3^{\mathcal{T}}$ be the random candidate support returned by Algorithm~\ref{alg:subspace_sampling}. Then
		\begin{align*}
			\mathbb{E}_{\mathcal{T}}
			\left[
			v_{I}(S_3^{\mathcal{T}})
			\right]
			\leq
			\left(
			1+\frac{d-t}{d}
			\right)
			v_I^{*}.
		\end{align*}
	\end{thm}
	\begin{proof}
		\begin{subequations}
			By the type-$1$ analogue of Theorem \ref{thm:ratio_support_reduction}, it suffices to prove the corresponding approximation guarantee for arbitrary support points $\bm{x}_i\in\Xi_i$ for all $i\in[k]$. 
			Let
			$
			\bm{c}
			\in
			\operatorname*{arg\,min}_{\bm{s}\in\mathbb{R}^d}
			\sum_{i\in[k]}\lambda_i\|\bm{s}-\bm{x}_i\|_1.
			$

			Since the $\ell_1$ objective is separable across coordinates, $\bm{c}$ can be chosen such that each coordinate $c_j$ is a weighted median of $\{x_{ij}\}_{i\in[k]}$ and, in particular, is equal to one of the observed coordinate values. Thus,
			$
			c_j\in U_j$ for all $j\in[d].$
			
			For a fixed realization of $\mathcal{T}$, define
			\begin{align*}
				D_{\mathcal{T}}(\bm{c})
				:=
				\min_{\bm{s}\in S_3^{\mathcal{T}}}
				\|\bm{s}-\bm{c}\|_1,
			\end{align*}
			and choose
			\begin{align*}
				\bm{s}_{\mathcal{T}}
				\in
				\operatorname*{arg\,min}_{\bm{s}\in S_3^{\mathcal{T}}}
				\|\bm{s}-\bm{c}\|_1.
			\end{align*}
			By the triangle inequality and $\sum_{i\in[k]}\lambda_i=1$, we obtain
			\begin{align}
				\min_{\bm{s}\in S_3^{\mathcal{T}}}\sum_{i\in[k]}\lambda_i\|\bm{s}-\bm{x}_i\|_1
				&\leq 
				\sum_{i\in[k]}
				\lambda_i
				\|\bm{s}_{\mathcal{T}}-\bm{x}_i\|_1 \leq
				\sum_{i\in[k]}
				\lambda_i
				\left(
				\|\bm{s}_{\mathcal{T}}-\bm{c}\|_1
				+
				\|\bm{c}-\bm{x}_i\|_1
				\right) \notag\\
				&=
				D_{\mathcal{T}}(\bm{c})+\sum_{i\in[k]}
				\lambda_i
				\|\bm{c}-\bm{x}_i\|_1 .
				\label{eq:W1_triangle_candidate}
			\end{align}
			
			We next bound the expected value of $D_{\mathcal{T}}(\bm{c})$. For each $i\in[k]$, define the mixed point $\tilde{\bm{x}}_i^{\mathcal{T}}\in\mathbb{R}^d$ coordinatewise by
			\begin{align*}
				\tilde{x}_{ij}^{\mathcal{T}}
				=
				\begin{cases}
					x_{ij}, & j\in\mathcal{T},\\
					c_j, & j\notin\mathcal{T}.
				\end{cases}
			\end{align*}
			The coordinates indexed by $\mathcal{T}$ are inherited from the single point $\bm{x}_i\in\hat{U}$, while each remaining coordinate $c_j$ belongs to $U_j$. Therefore,
			$
			\tilde{\bm{x}}_i^{\mathcal{T}}
			\in
			S_3^{\mathcal{T}}.
			$
			It follows that
			\begin{align}
				D_{\mathcal{T}}(\bm{c})
				\leq
				\left\|
				\tilde{\bm{x}}_i^{\mathcal{T}}-\bm{c}
				\right\|_1
				=
				\sum_{j\in\mathcal{T}}
				|x_{ij}-c_j|,
				\qquad
				\forall i\in[k].
				\label{eq:W1_candidate_distance}
			\end{align}
			
			Since $\mathcal{T}$ is sampled uniformly from all subsets of $[d]$ of cardinality $d-t$, every coordinate belongs to $\mathcal{T}$ with probability $(d-t)/d$. Hence, for every $i\in[k]$,
			\begin{align}
				\mathbb{E}_{\mathcal{T}}
				\left[
				\left\|
				\tilde{\bm{x}}_i^{\mathcal{T}}-\bm{c}
				\right\|_1
				\right]
				&=
				\mathbb{E}_{\mathcal{T}}
				\left[
				\sum_{j\in\mathcal{T}}
				|x_{ij}-c_j|
				\right] =
				\sum_{j\in[d]}
				\Pr[j\in\mathcal{T}]
				|x_{ij}-c_j| =
				\frac{d-t}{d}
				\sum_{j\in[d]}
				|x_{ij}-c_j| \notag\\
				&=
				\frac{d-t}{d}
				\|\bm{x}_i-\bm{c}\|_1.
				\label{eq:expected_l1_norm}
			\end{align}
			Taking expectations in \eqref{eq:W1_candidate_distance}, multiplying the resulting inequality by $\lambda_i$, and summing over $i\in[k]$ yield
			\begin{align}
				\mathbb{E}_{\mathcal{T}}
				\left[
				D_{\mathcal{T}}(\bm{c})
				\right]
				&=
				\sum_{i\in[k]}
				\lambda_i
				\mathbb{E}_{\mathcal{T}}
				\left[
				D_{\mathcal{T}}(\bm{c})
				\right] \leq
				\sum_{i\in[k]}
				\lambda_i
				\mathbb{E}_{\mathcal{T}}
				\left[
				\left\|
				\tilde{\bm{x}}_i^{\mathcal{T}}-\bm{c}
				\right\|_1
				\right] =
				\frac{d-t}{d}
				\sum_{i\in[k]}
				\lambda_i
				\|\bm{x}_i-\bm{c}\|_1 
				\label{eq:expected_distance_candidate_set}
			\end{align}
			
			Taking expectations in \eqref{eq:W1_triangle_candidate} and applying \eqref{eq:expected_distance_candidate_set}, we obtain
			\begin{align*}
				\mathbb{E}_{\mathcal{T}}
				\left[
				\min_{\bm{s}\in S_3^{\mathcal{T}}}
				\sum_{i\in[k]}\lambda_i\|\bm{s}-\bm{x}_i\|_1
				\right]
				&\leq
				\sum_{i\in[k]} \lambda_i \|\bm{x}_i-\bm{c}\|_1 
				+
				\mathbb{E}_{\mathcal{T}}
				\left[
				D_{\mathcal{T}}(\bm{c})
				\right]\leq
				\left(
				1+\frac{d-t}{d}
				\right)
				\sum_{i\in[k]} \lambda_i \|\bm{x}_i-\bm{c}\|_1 .
			\end{align*}
			Since the tuple $(\bm{x}_1,\ldots,\bm{x}_k)$ was arbitrary, the result follows from Theorem \ref{thm:ratio_support_reduction}.
		\end{subequations}
	\end{proof}
	
	\paragraph{Deterministic counterpart: subspace enumeration.}
	Enumerating all subsets $\mathcal T\subseteq[d]$ with $|\mathcal T|=d-t$ gives
	$S_3^t
	:=
	\bigcup_{\substack{\mathcal{T}\subseteq[d]\\|\mathcal{T}|=d-t}}
	S_3^{\mathcal{T}}.
	$
	Its cardinality is bounded by
	\begin{align*}
		|S_3^t|
		\leq
		\binom{d}{d-t}(nk)^{t+1}
		=
		\binom{d}{t}(nk)^{t+1}
		=O\left(d^t(nk)^{t+1}\right)
	\end{align*}
	for fixed $t$. Then we have the following determinsitic performance guarantee.
	
	\begin{algorithm}[htbp]
		\caption{Subspace-enumeration algorithm for the type-$1$ Wasserstein barycenter problem}
		\label{alg:subspace_enumeration}
		\begin{algorithmic}[1]
			\State \textbf{Input}: Probability measures $\{\mathbb{P}_i\}_{i\in[k]}$ with supports $\{\Xi_i\}_{i\in[k]}$ and an integer $t\in\{0,\ldots,d\}$.
			\State Construct
			$
			S_3^t
			:=
			\bigcup_{\substack{\mathcal{T}\subseteq[d]\\|\mathcal{T}|=d-t}}
			S_3^{\mathcal{T}}.
			$
			\State Solve the type-$1$ restricted MOT problem over $S_3^t$ and recover an optimal primal solution $\bm{\Pi}^*(S_3^t)$.
			\State \textbf{Output}: Construct and return the approximate barycenter $\mathbb{P}(S_3^t)$ using \Cref{lem_distr}.
		\end{algorithmic}
	\end{algorithm}
	
	\begin{thm}\label{thm:W1_enumeration}
		Let $S_3^t$ be the candidate support returned by Algorithm~\ref{alg:subspace_enumeration}. Then
		\begin{align*}
			v_{I}(S_3^t)
			\leq
			\left(
			1+\frac{d-t}{d}
			\right)
			v_{I}^{*}.
		\end{align*}
	\end{thm}
	\begin{proof}
		\begin{subequations}
			As in the proof of Theorem \ref{thm:W1_sampling}, fix arbitrary support points $\bm{x}_i\in\Xi_i$, $i\in[k]$, and let
			\begin{align*}
				\bm{c}
				\in
				\operatorname*{arg\,min}_{\bm{s}\in\mathbb{R}^d}
				\sum_{i\in[k]}
				\lambda_i\|\bm{s}-\bm{x}_i\|_1
			\end{align*}
			be chosen such that $c_j$ is an observed coordinate value for every $j\in[d]$. 
			For each $i\in[k]$, choose a subset
			\begin{align*}
				\mathcal{T}_i
				\in
				\operatorname*{arg\,min}_{\substack{\mathcal{T}\subseteq[d]\\|\mathcal{T}|=d-t}}
				\sum_{j\in\mathcal{T}}
				|x_{ij}-c_j|.
			\end{align*}
			Thus, $\mathcal{T}_i$ consists of the $d-t$ coordinates having the smallest values of $|x_{ij}-c_j|$. Therefore,
			\begin{align}
				\sum_{j\in\mathcal{T}_i}
				|x_{ij}-c_j|
				\leq
				\frac{d-t}{d}
				\sum_{j\in[d]}
				|x_{ij}-c_j|
				=
				\frac{d-t}{d}
				\|\bm{x}_i-\bm{c}\|_1.
				\label{eq:deterministic_coordinate_average}
			\end{align}
			
			Define $\tilde{\bm{x}}_i\in\mathbb{R}^d$ by
			\begin{align*}
				\tilde{x}_{ij}
				=
				\begin{cases}
					x_{ij}, & j\in\mathcal{T}_i,\\
					c_j, & j\notin\mathcal{T}_i.
				\end{cases}
			\end{align*}
			By construction, we have
			$
			\tilde{\bm{x}}_i
			\in
			S_3^{\mathcal{T}_i}
			\subseteq
			S_3^t.
			$
			Therefore, for every $i\in[k]$,
			\begin{align}
				\min_{\bm{s}\in S_3^t}
				\|\bm{s}-\bm{c}\|_1
				&\leq
				\|\tilde{\bm{x}}_i-\bm{c}\|_1 =
				\sum_{j\in\mathcal{T}_i}
				|x_{ij}-c_j| \leq
				\frac{d-t}{d}
				\|\bm{x}_i-\bm{c}\|_1.
				\label{eq:deterministic_candidate_distance}
			\end{align}
			Multiplying \eqref{eq:deterministic_candidate_distance} by $\lambda_i$ and summing over $i\in[k]$ gives
			\begin{align}
				\min_{\bm{s}\in S_3^t}
				\|\bm{s}-\bm{c}\|_1
				\leq
				\frac{d-t}{d}
				\sum_{i\in[k]}
				\lambda_i
				\|\bm{x}_i-\bm{c}\|_1
				=
				\frac{d-t}{d}
				\sum_{i\in[k]} \lambda_i \|\bm{x}_i-\bm{c}\|_1 .
				\label{eq:deterministic_set_distance}
			\end{align}
			
			Let
			$
			\bm{s}^*
			\in
			\operatorname*{arg\,min}_{\bm{s}\in S_3^t}
			\|\bm{s}-\bm{c}\|_1.
			$
			Using the triangle inequality and \eqref{eq:deterministic_set_distance}, we obtain
			\begin{align*}
				\min_{\bm{s}\in S_3^t}\sum_{i\in[k]}\lambda_i\|\bm{s}-\bm{x}_i\|_1
				&\leq
				\sum_{i\in[k]} \lambda_i \|\bm{x}_i-\bm{s}^*\|_1 \leq
				\sum_{i\in[k]} \lambda_i \|\bm{x}_i-\bm{c}\|_1 
				+
				\|\bm{s}^*-\bm{c}\|_1\\
				&=
				\sum_{i\in[k]} \lambda_i \|\bm{x}_i-\bm{c}\|_1 
				+
				\min_{\bm{s}\in S_3^t}
				\|\bm{s}-\bm{c}\|_1\leq
				\left(
				1+\frac{d-t}{d}
				\right)
				\sum_{i\in[k]} \lambda_i \|\bm{x}_i-\bm{c}\|_1 .
			\end{align*}
			Since this inequality holds for every tuple $(\bm{x}_1,\ldots,\bm{x}_k)$, the result follows from Theorem \ref{thm:ratio_support_reduction}.
		\end{subequations}
	\end{proof}

	\color{black}
	
	\section{Conclusion and future work}\label{section:conclusion}
	We developed a candidate-support reduction framework for discrete Wasserstein barycenters. For type-$2$ barycenters with general weights, the framework yields a randomized approximation scheme with an
	expected guarantee and a deterministic PTAS, both of which generalize the previously known factor-$2$ method. Under equal weights, a sampling
	procedure without replacement leads to a sharper approximation ratio. The same candidate-support principle yields a $(1+1/t)$ support-reduction
	guarantee for sparse type-$2$ barycenters and both exact and approximate
	coordinate-based constructions for type-$1$ barycenters. The numerical results show the efficiency of the proposed approximation algorithms. Important directions for future work include extending the approximation framework
	to type-$p$ costs with $p\notin\{1,2\}$, developing more scalable methods
	for solving the restricted MOT problems, and identifying other
	optimization problems for which tuplewise candidate-support reduction
	provides provable approximation guarantees.

	\section*{Acknowledgments}
	This research was supported in part by National Science Foundation grant 2246414 and Office of Naval Research grant N00014-24-1-2066. The second author thanks Zedong Wang and Shuai Li for helpful discussions.

	\bibliographystyle{agsm}
	\bibliography{refs.bib}
	
\end{document}